\documentclass[12pt,letterpaper]{amsart}
\usepackage[latin1]{inputenc}

\makeatletter

\theoremstyle{plain}    
\newtheorem{thm}{Theorem}[section]
\numberwithin{equation}{section} 
\numberwithin{figure}{section} 
\theoremstyle{plain}    
\newtheorem{cor}[thm]{Corollary} 
\theoremstyle{plain}    
\newtheorem{lem}[thm]{Lemma} 
\theoremstyle{plain}    
\newtheorem{prop}[thm]{Proposition} 
\theoremstyle{plain}    
\newtheorem*{prop*}{Proposition} 
\theoremstyle{definition}
\newtheorem{defn}[thm]{Definition}
\theoremstyle{remark}
\newtheorem{rem}[thm]{Remark}
\theoremstyle{remark}    
\newtheorem{claim}[thm]{Claim}
\theoremstyle{remark}    
\newtheorem*{acknowledgement*}{Acknowledgement} 

\usepackage{pslatex}
\usepackage{amssymb}

\def\btw{\between}
\def\id{\operatorname{id}}
\def\Tr{\operatorname{Tr}}
\def\Ad{\operatorname{Ad}}

\makeatother

\begin{document}

\title{Microstates free entropy and cost of equivalence relations.}

\begin{abstract}
We define an analog of Voiculescu's free entropy for \( n \)-tuples of unitaries
\( u_{1},\dots ,u_{n} \) in a tracial von Neumann algebra \( M \), normalizing
a unital subalgebra \( L^{\infty }[0,1]=B\subset M \). Using this quantity,
we define the free dimension \( \delta _{0}(u_{1},\dots ,u_{n}\btw B) \). This
number depends on \( u_{1},\dots ,u_{n} \) only up ``orbit equivalence''
over \( B \). In particular, if \( R \) is an measurable equivalence relation
on \( [0,1] \) generated by \( n \) automorphisms \( \alpha _{1},\dots ,\alpha _{n} \),
let \( u_{1},\dots ,u_{n} \) be the unitaries implementing \( \alpha _{1},\dots ,\alpha _{n} \)
in the Feldman-Moore crossed product algebra \( M=W^{*}([0,1],R)\supset B=L^{\infty }[0,1] \).
In this way, we obtain an invariant \( \delta (R)=\delta _{0}(u_{1},\dots ,u_{n}\btw B) \)
of the equivalence relation \( R \). If \( R \) is treeable, \( \delta (R) \)
coincides with the cost \( C(R) \) of \( R \) in the sense of Gaboriau. For
a general equivalence relation \( R \) posessing a finite graphing, \( \delta (R)\leq C(R) \).
Using the notion of free dimension, we define an dynamical entropy invariant
for an automorphism of a measurable equivalence relation (or more generally
of an \( r \)-discrete measure groupoid), and give examples.
\end{abstract}

\author{Dimitri Shlyakhtenko}

\date{\today}

\address{Department of Mathematics, UCLA, Los Angeles, CA 90095}

\email{shlyakht@math.ucla.edu }

\thanks{Research supported by an NSF postdoctoral fellowship.}

\maketitle

\section{Introduction.}

This is our second paper investigating the connections between the notion of
cost of equivalence relations introduced by Gaboriau in \cite{gaboriau:cost:cr},
\cite{gaboriau:cost} and Voiculescu's free entropy and free dimension theory
\cite{dvv:entropy1}, \cite{dvv:entropy2}, \cite{dvv:entropy3}, \cite{dvv:entropy4},
\cite{dvv:entropy5}, \cite{dvv:entropy6} and \cite{dvv:improvedrandom}. While
in our first paper \cite{shlyakht:cost} we used the non-commutative Hilbert
transform approach to free entropy (the so-called microstates-free approach),
this paper is concerned with the microstates approach. 

If \( M \) is a tracial von Neumann algebra, and \( L^{\infty }[0,1]\cong B\subset M \)
is a unital \( W^{*} \)-subalgebra, we associate to each \( n \)-tuple of
unitaries \( u_{1},\dots ,u_{n} \) in the normalizer of \( B \) its entropy
with respect to \( B \), \( \chi (u_{1},\dots ,u_{n}\btw B) \). In the first
approximation, \( \chi  \) measures the extent to which \( u_{1},\dots ,u_{n} \)
are free with amalgamation over \( B \). We caution the reader that \( \chi (u_{1},\dots ,u_{n}\btw B) \)
is not the entropy of \( u_{1},\dots ,u_{n} \) relative to \( B \). Indeed,
such a relative entropy must measure freeness between \( u_{1},\dots ,u_{n} \)
and \( B \). In our case, \( u_{j}Bu_{j}^{*}=B \), since \( u_{j} \) are
assumed to normalize \( B \), and hence \( u_{1},\dots ,u_{n} \) cannot be
free from \( B \).

Using \( \chi (\cdots \btw B) \), we define in the spirit of Voiculescu's definition
of free dimension the quantity
\[
\delta _{0,\kappa }^{\omega }(u_{1},\dots ,u_{n}\btw B),\]
which we call the free dimension of \( u_{1},\dots ,u_{n} \) with respect to
\( B \). We show that \( \delta _{0,\kappa }^{\omega } \) depends on \( u_{1},\dots ,u_{n} \)
only up to \char`\"{}orbit equivalence'' over \( B \). Furthermore, \( \delta _{0,\kappa }^{\omega }(u_{1},\dots ,u_{n},v_{1}\dot{,}\dots ,v_{m}\btw B)=\delta _{0,\kappa }^{\omega }(u_{1},\dots ,u_{n}\btw B)+\delta _{0,\kappa }^{\omega }(v_{1},\dots ,v_{m}\btw B) \)
if \( (u_{1},\dots ,u_{n}) \) and \( (v_{1},\dots ,v_{m}) \) are free with
amalgamation over \( B \). We explicitely compute \( \delta _{0,\kappa }^{\omega }(u\btw B) \)
in the case that \( u \) is the implementing unitary for a free measure-preserving
action of a cyclic group on \( B \). 

If \( R \) is a measurable measure-preserving equivalence relation on \( [0,1] \),
Feldman and Moore associated to it a von Neumann algebra \( W^{*}([0,1],R)=M \)
(see \cite{feldman-moore}). In the case that \( R \) can be generated by \( n \)
automorphisms \( \alpha _{1},\dots ,\alpha _{n} \) (i.e., has a ``graphing''
by \( \alpha _{1},\dots ,\alpha _{n} \)), \( W^{*}([0,1],R) \) is generated
by \( B=L^{\infty }[0,1] \) and unitaries \( u_{1},\dots ,u_{n} \) implementing
the automorphisms \( \alpha _{1},\dots ,\alpha _{n} \). Two choices of graphings
\( \alpha _{1},\dots ,\alpha _{n} \) and \( \beta _{1},\dots ,\beta _{m} \)
give rise to two families of unitaries \( u_{1},\dots ,u_{n} \) and \( v_{1},\dots ,v_{m} \),
which are orbit-equivalent over \( B \). It follows that
\[
\delta _{0,\kappa }^{\omega }(B\subset M)=\delta ^{\omega }_{0,\kappa }(u_{1},\dots ,u_{n}\btw B)=\delta _{0,\kappa }^{\omega }(v_{1},\dots ,v_{m}\btw B)\]
is independent of the choice of the graphing \( \alpha _{1},\dots ,\alpha _{n} \),
and is an invariant of the pair \( B\subset M \). This invariant satisfies
\[
\delta _{0,\kappa }^{\omega }(B\subset M)=\delta _{0,\kappa }^{\omega }(B\subset M_{1})+\delta _{0,\kappa }^{\omega }(B\subset M_{2})\]
if \( M=W^{*}(M_{1},M_{2},B) \) and \( M_{1},M_{2}\subset M \) are free with
amalgamation over \( B \). 

In particular, for \( B=L^{\infty }[0,1]\subset M=W^{*}([0,1],R) \), the number
\( \delta _{0,\kappa }^{\omega }(B\subset M) \) is an invariant of the equivalence
relation \( R \). Let us write \( \delta (R) \) for its value. Gaboriau recently
introduced another invariant of an equivalence relation, which he calls the
cost \( C(R) \) (see \cite{gaboriau:cost:cr}, \cite{gaboriau:cost}). If \( R \)
is generated by two sub-equivalence relations \( R_{1} \) and \( R_{2} \),
such that \( R_{1} \) and \( R_{2} \) are free inside \( R \), then he proved
that \( C(R)=C(R_{1})+C(R_{2}) \). We show that if \( R \) is an equivalence
relation generated by a single automorphism, then \( \delta (R)=C(R) \). This
means that if \( R \) is an arbitrary treeable equivalence relation (i.e.,
\( R \) is generated by a family of singly-generated subrelations \( R_{i} \)
with \( R_{i} \) free), then \( C(R)=\delta (R) \). In general, we have \( C(R)\geq \delta (R) \).
It is possible that in fact one has \( C(R)=\delta (R) \); however, this would
in particular imply that an arbitrary von Neumann algebra having a Cartan subalgebra
can be embedded into an ultrapower of the hyperfinite II\( _{1} \) factor.

We mention that there is a similarity between properties of microstates and
microstates-free entropies. Therefore, one may expect that many properties of
microstates-free free entropy and free dimension with respect to \( B \) \cite{shlyakht:cost}
should have analogs for the microstates quantities considered in the present
paper. In particular, consider the following proposition from \cite{shlyakht:cost};
here \( \delta ^{*} \) refers to the microstates-free free dimension:

\begin{prop*}
Let \( \alpha  \) be a free measure-preserving action of a group \( G \) on
\( [0,1] \). Assume that \( g_{1},\dots ,g_{n}\in G \) generate \( G \).
Let \( R \) be the equivalence relation induced by this action, and let \( u_{1},\dots ,u_{n} \)
be unitaries in \( W^{*}([0,1],R)\supset B=L^{\infty }[0,1] \) corresponding
to \( \alpha _{g_{1}},\dots ,\alpha _{g_{n}} \). Let \( v_{1},\dots ,v_{n} \)
be unitaries in the group von Neumann algebra of \( G \), corresponding to
the generators \( g_{1},\dots ,g_{n} \). Then
\[
\delta ^{*}(u_{1},\dots ,u_{n}\btw B)=\delta ^{*}(v_{1},\dots ,v_{n})\]
and in particular depends only on \( g_{1},\dots ,g_{n}\in G \).
\end{prop*}
If this proposition were to hold for \( \delta ^{\omega }_{0,\kappa } \) instead
of \( \delta ^{*} \), we would obtain that \( \delta (R)=\delta (G) \), where
\( G \) is any finitely-generated group, \( R \) is a measurable equivalence
relation induced by an \emph{arbitrary} free action of \( G \) on a finite
measure space, and \( \delta (G) \) refers to the free dimension of \( G \)
introduced by Voiculescu in \cite{dvv:improvedrandom}. In particular, this
would give \( \delta (G)\leq C(G) \) for any finitely-generated group \( G \).

The free dimension \( \delta (R) \) measures the ``size'' of the equivalence
relation \( R \). Using this, we define an entropy-like invariant for a dynamical
system involving automorphisms of equivalence relations. For a free shift of
multiplicity \( n \), this invariant is \( n \).

\begin{acknowledgement*}
This work was carried out while visiting Centre Émile Borel, Institut Henri
Poicaré, Paris, France, to which I am grateful for the friendly and encouraging
atmosphere. I would like to especially thank the organizers of the Free Probability
and Operator Spaces program at IHP, Professors P. Biane, G. Pisier and D. Voiculescu,
for a very stimulating semester.
\end{acknowledgement*}

\section{Preliminaries and Notation.\label{sec:notation}}

\subsection{Basic notation.}

We denote by \( M_{N\times N} \) the algebra of complex \( N\times N \) matrices,
and by \( \Delta _{N} \) its subalgebra consisting of diagonal matrices. 

Note that \( \Delta _{N}\subset L^{\infty }[0,1] \) as the algebra of functions,
which are piece-wise constant on the intervals \( [\frac{k}{N},\frac{k+1}{N}] \),
\( 0\leq k<N \). We denote by \( S_{N} \) the symmetric group of permutations
of size \( N \); \( S_{N} \) acts on \( \Delta _{N} \) in the obvious way.
We denote by \( U(N) \) the unitary group of \( M_{N\times N} \).

We denote by \( \Tr  \) the usual matrix trace on \( M_{N\times N} \); \( \Tr (I)=N \),
where \( I \) denotes the identity matrix. 

Although it should always be clear from the context, we try to adhere to the
following general notational rule: elements of \( M_{N\times N} \) will be
denoted by capital letters (\( U,V \), etc.), while elements of abstract von
Neumann algebras will be denoted by lower-case letters (\( u,v \), etc.).

\subsection{Operator-valued distributions.}

We recall some standard notions from free probability theory (see \cite{DVV:book},
\cite{dvv:amalg} for more details). Let \( M \) be a von Neumann algebra,
\( \tau  \) be a faithful state on \( M \) and \( B \) be a unital von Neumann
subalgebra. Then there always exists a conditional expectation \( E=E_{B}:M\to B \),
determined by:
\[
E(bmb')=bE(m)b',\qquad b,b\in B,\quad m\in M\]

\[
\tau (bm)=\tau (bE(m)),\qquad b\in B,\quad m\in M.\]
If \( u_{1},\dots ,u_{n}\in M \) is a family of elements, we refer to each
expression
\[
E_{B}(b_{0}u_{i_{1}}b_{1}\cdots u_{i_{n}}b_{n})\]
as a \( B \)-valued moment of \( (u_{1},\dots ,u_{n}) \). The moments define
a linear map \( \mu _{(u_{1},\dots ,u_{n})} \) from the algebra \( B[t_{1},\dots ,t_{n}] \)
of non-commutative polynomials with coefficients from \( B \) on \( n \) non-commuting
indeterminates to \( B \) by
\[
\mu _{(u_{1},\dots ,u_{n})}(b_{0}t_{i_{1}}b_{1}\cdots t_{i_{n}}b_{n})=E_{B}(b_{0}u_{i_{1}}b_{1}\cdots u_{i_{n}}b_{n}).\]
If the variables \( u_{1},\dots ,u_{n} \) are not self-adjoint, we refer to
the distribution of the family \( (u_{1},u_{1}^{*},\dots ,u_{n},u_{n}^{*}) \)
as the \( * \)-distribution of \( u_{1},\dots ,u_{n} \).

The \( B \)-valued \( * \)-distribution of \( (u_{1},\dots ,u_{n}) \) determines
(up to isomorphism) the pair \( B\subset W^{*}(B,u_{1},\dots ,u_{n}) \) (here
for a set \( S \), \( W^{*}(S) \) denotes the von Neumann algebra generated
by \( S \)).

In the case \( B=\mathbb {C} \) we speak of a distribution (or \( * \)-distribution)
of a family. Note that the knowledge of the \( B \)-valued distribution of
a family \( (u_{1},\dots ,u_{n}) \) is equivalent to knowledge of the \( \mathbb {C} \)-valued
distribution of \( (u_{1},\dots ,u_{n},b_{1},\dots ,b_{n},\dots ) \) where
\( b_{1},b_{2},\dots  \) are some generators of \( B \).

As an example, consider the algebra \( M=M_{dN\times dN}=M_{d\times d}\otimes M_{N\times N} \)
of \( dN\times dN \) matrices, and in it the subalgebra \( M_{N\times N}\cong B=1\otimes M_{N\times N} \)
. Then each element of \( M \) can be written as an \( N\times N \) block
matrix, with blocks of size \( d\times d \). The knowledge of the \( B \)-valued
distribution of some family of matrices \( U_{1},\dots ,U_{n}\in M \) is equivalent
to the knowledge of the joint (scalar) distribution of their constituent blocks.

\subsection{Freeness with amalgamation.}

Let \( M \) be a von Neumann algebra with a faithful trace \( \tau  \), and
\( B \) be a von Neumann subalgebra. Denote by \( E \) the canonical \( B \)-valued
conditional expectation onto \( B \). Let \( M_{i}\subset M \) be subalgebras
containing \( B \). Then \( M_{i} \) are free with amalgamation over \( B \)
if
\[
E(m_{1}\dots m_{n})=0\]
whenever \( m_{j}\in M_{i(j)} \), \( i(1)\neq i(2) \), \( \ldots  \), \( i(n-1)\neq i(n) \),
and \( E(m_{j})=0 \).

We say that sets \( X_{1},\dots ,X_{n}\subset M \) are \( * \)-free over \( B \)
if the algebras \( M_{i}=W^{*}(X_{i},B) \) are free with amalgamation over
\( B \).

If families \( X_{1}=(u_{1},\dots ,u_{n}) \), \( X_{2} \), \( \ldots  \),
\( X_{n} \) are \( * \)-free over \( B \), then the joint \( B \)-valued
\( * \)-distribution of \( \sqcup X_{j} \) is completely determined by the
\( B \)-valued \( * \)-distributions of each family \( X_{j} \).

\subsection{Independence.}

Let \( M \) be a von Neumann algebra with a faithful trace \( \tau  \), and
\( B \) be a von Neumann subalgebra. Denote by \( E \) the canonical \( B \)-valued
conditional expectation onto \( B \). Let \( A\subset M \) be a subalgebra.
Then \( A \) are independent from \( B \), if:
\[
[a,b]=0,\quad \tau (ab)=\tau (a)\tau (b),\qquad \forall a\in A,b\in B.\]
If \( X=(u_{1},\dots ,u_{n})\in M \) is a family of variables, then we say
that \( X \) is independent from \( B \), if the algebra \( A=W^{*}(u_{1},\dots ,u_{n}) \)
is independent from \( B \). Note that if a family \( X \) is independent
from \( B \), then then its \( B \)-valued \( * \)-distribution is determined
completely by its scalar-valued \( * \)-distribution.

\subsection{Normalizer \protect\( \mathcal{N}(B)\protect \). }

If \( B\subset M \) is a diffuse commutative von Neumann subalgebra, we denote
by \( \mathcal{N}(B) \) the set
\[
\mathcal{N}(B)=\{u\in M\textrm{ unitary}:uBu^{*}\subset B\}.\]
Unitaries in \( \mathcal{N}(B) \) are said to \emph{normalize} \( B \).

\subsection{Preliminaries on \protect\( |\cdot |_{\epsilon }\protect \).}

We will be concerned with approximating \( L^{\infty }[0,1] \)-valued distributions
of non-commutative random variables with distributions of matrices. It will
be useful to introduce the following quantity:

\begin{defn}
Let \( f\in L^{\infty }[0,1] \) and let \( \epsilon >0 \). Then 
\[
|f|_{\epsilon }=\inf _{X\subset [0,1],\mu (X)\geq 1-\epsilon }\sup _{\xi \in X}|f(\xi )|.\]

\end{defn}
\begin{lem}
One has \( |\alpha f|_{\epsilon }=|\alpha ||f|_{\epsilon } \), \( |f|_{\epsilon }\leq |f|_{\delta } \)
if \( \epsilon \geq \delta  \); \( |f+g|_{\epsilon +\epsilon '}\leq |f|_{\epsilon }+|g|_{\epsilon '} \)
and \( |fg|_{\epsilon +\epsilon '}\leq |f|_{\epsilon }|g|_{\epsilon '} \).
\end{lem}
Note that the family \( |\cdot |_{\epsilon } \) induces a topology \( \tau  \)
on \( L^{\infty }[0,1] \): a sequence of functions \( \{f_{n}\} \) converges
to a function \( g \) iff \( |f_{n}-g|_{\epsilon }\to 0 \) for all \( \epsilon >0 \).

\begin{lem}
\label{lemma:sametop}The topology \( \tau  \) coninsides with the topology
of strong convergence in \( L^{2}[0,1] \) on \( \Vert \cdot \Vert _{\infty } \)-bounded
subsets of \( L^{\infty }[0,1] \).
\end{lem}
\begin{rem}
Every function \( d\in L^{\infty }[0,1] \) can be approximated in \( |\cdot |_{\epsilon } \)
by functions from \( \Delta _{N} \) for \( N \) sufficiently large. Indeed,
it is sufficient to show that any step-function \( d \), which is constant
on subsets \( X_{1},\dots ,X_{n} \) of \( [0,1] \), can be approximated in
this way. But this is equivalent to showing that there exists \( N \) sufficiently
large, and disjoint subsets \( S_{1},\dots ,S_{n} \) of \( \{0,\dots ,N-1\} \),
so that if we set \( Y_{j}=\sqcup _{k\in S_{j}}[k/N,(k+1)/N] \), one has that
\( \bigcup _{j}\left( (X_{j}\setminus Y_{j})\cup (Y_{j}\setminus X_{j})\right)  \)
has measure less than \( \epsilon  \). In fact, \( d_{n} \) can be chosen
so that \( \Vert d_{j}\Vert _{\infty }\leq \Vert d\Vert _{\infty } \). (This
remark can also be seen from strong density of \( \cup _{N}\Delta _{N} \) in
the unit ball of \( L^{\infty }[0,1] \)).
\end{rem}

\subsection{Some approximation lemmas.}

\begin{lem}
Let \( \sigma :[0,1]\to [0,1] \) be a measure-preserving Borel isomorphism.
Then, given \( \epsilon ,\delta >0 \) and \( N_{0}>0 \) and \( d\in L^{\infty }[0,1] \),
there exists \( N>N_{0} \), and a permutation \( \Sigma \in S_{N} \), so that
\[
|\sigma (d)-\Sigma (d)|_{\epsilon }<\delta .\]

\end{lem}
\begin{proof}
For two partitions \( P,Q \) we say that \( |P-Q|<\epsilon  \), if there exists
a set \( Y\subset [0,1] \) of measure \( \lambda (Y)\geq 1-\epsilon  \), and
such that \( P_{i}\cap Y \) and \( Q_{j}\cap Y \) are either distinct, or
coinside, for all \( i,j \). 

Denote by \( [M] \) the partition \( \{[0,\frac{1}{M}],[\frac{1}{M},\frac{2}{M}],\dots ,[\frac{M-1}{M},1]\} \)
of \( [0,1] \). 

Note that for any partition \( P \), there exists an \( M>N_{0} \), such that
\( |P-[M]|<\epsilon  \).

One can assume, by replacing \( \delta  \) with \( \lambda \delta  \), \( \lambda >0 \),
that \( \Vert d\Vert _{\infty }\leq 1 \); one can also assume that \( 1/\delta  \)
is an integer. 

Let \( P \) be the partition of \( [0,1] \) given by
\[
P_{j}=d^{-1}([j\delta /8,(j+1)\delta /8]),\quad -8/\delta \leq j<8/\delta .\]
 It follows that \( \sup _{\xi ,\zeta \in P_{j}}|d(\xi )-d(\zeta )|\leq \delta /4 \). 

For \( M>N_{0} \) sufficiently large, \( |P-[M]|<\epsilon /4 \) and \( |\sigma (P)-[M]|<\epsilon /4 \);
hence there exists a step-function \( d'\in \Delta _{M} \) for which \( |d-d'|_{\epsilon /4}<\delta /4 \);
we may also require there exists a step-function \( \delta ''\in \Delta _{M} \),
for which \( |\sigma (d)-d''|_{\epsilon /4}<\delta /4 \). There is a permutation
\( \Sigma \in S_{M} \), so that \( |d''-\Sigma (d')|_{\epsilon /2}<\delta /2 \);
this is because the measure of \( (d'')^{-1}([j\delta /8,(j+1)\delta /8])\cap Y \)
and \( (d')^{-1}([j\delta /8,(j+1)\delta /8)\cap Y \) is the same, for a set
\( Y\subset [0,1] \) of measure \( \geq 1-\epsilon /4 \). It follows that
\( |\sigma (d)-\Sigma (d')|_{3\epsilon /4}\leq |\sigma (d)-d''|_{\epsilon /4}+|d''-\Sigma (d')|_{\epsilon /2}<3\delta /4 \).
Finally, we get that
\[
|\sigma (d)-\Sigma (d)|_{\epsilon }\leq |\sigma (d)-\Sigma (d')|_{3\epsilon /4}+|\Sigma (d')-\Sigma (d)|_{\epsilon /4}<3\delta /4+|d'-d|_{\epsilon /4}=\delta .\]

\end{proof}
\begin{cor}
\label{corr:existanceofpermutation}Let \( \sigma  \) be a measure preserving
automorphism of \( [0,1] \). Fix \( N>0 \), \( \epsilon >0 \), \( \delta >0 \),
\( l>0 \) and \( d_{0},d_{1},\dots ,d_{l}\in L^{\infty }[0,1] \). Then there
exists an \( M>N \) and a permutation \( \Sigma \in S_{M} \), for which
\[
|d_{0}\sigma ^{g(1)}(d_{1}\sigma ^{g(1)}(\dots \sigma ^{g(k)}(d_{k})\dots )))-d_{0}\Sigma ^{g(1)}(d_{1}\Sigma ^{g(1)}(\dots \Sigma ^{g(k)}(d_{k})\dots )))|_{\epsilon }<\delta \]
 for all \( 1\leq k\leq l \) and \( g:\{1,\dots ,k\}\to \{\pm 1\} \).
\end{cor}

\subsection{Some freeness lemmas.}

\begin{lem}
Let \( B \) be an tracial von Neumann algebra. Let \( M=M_{N\times N}\otimes B \)
be the von Neumann algebra of \( B \)-valued \( N\times N \) matrices, with
the obvious trace. Let \( D\subset B \) be a subalgebra. Let \( m_{1},\dots ,m_{n}\in M \)
be elements, so that each \( m_{k} \) is a matrix
\[
m_{k}=(b_{ij}^{(k)})_{1\leq i,j\leq N},\quad 1\leq k\leq n.\]
Then \( m_{1},\dots ,m_{n} \) are \( * \)-free with amalgamation over \( M_{N\times N}\otimes D\subset M \)
in \( M \) if and only if the families \( F_{1}=(b_{ij}^{(1)}:1\leq i,j\leq N) \),
\( F_{2}=(b_{ij}^{(2)}:1\leq i,j\leq N) \), \( \ldots  \), \( F_{n}=(b_{ij}^{(n)}:1\leq i,j\leq N) \)
are \( * \)-free with amalgamation over \( D \). In particular, if \( D=\mathbb {C} \),
then \( m_{k} \) are \( * \)-free with amalgamation over \( M_{N\times N}\otimes 1 \)
iff the families \( F_{j} \), \( 1\leq j\leq n \) of their entries are \( * \)-free.
\end{lem}
The proof is a straightforward application of the freeness condition and matrix
multiplication.

\begin{lem}
\label{lemma:Ufreediag}Let \( B \) be an tracial von Neumann algebra. Let
\( M=M_{N\times N}\otimes B \) be the von Neumann algebra of \( B \)-valued
\( N\times N \) matrices, with the obvious trace. Let \( D\subset B \) be
a subalgebra. Let \( u_{1},\dots ,u_{N}\in B \) be unitaries, so that: \( u_{1},\dots ,u_{N} \)
are \( * \) -free with amalgamation over \( D \); each \( u_{j} \) is a Haar
unitary (i.e., \( \tau (u_{j}^{k})=0 \) unless \( k=0 \)); and \( u_{1},\dots ,u_{N} \)
are independent from \( D \). Let \( C\subset B \) be another subalgebra,
so that \( u_{1},\dots ,u_{n} \) are \( * \)-free from \( C \) with amalgamation
over \( D \). Consider in \( M \) the matrix
\[
U=\left( \begin{array}{ccc}
u_{1} & \cdots  & 0\\
\vdots  & \ddots  & \vdots \\
0 & \cdots  & u_{n}
\end{array}\right) .\]
Then \( U \) is \( * \)-free from \( M_{N\times N}\otimes C \) with amalgamation
over the algebra \( \Delta (D) \) of diagonal matrices with entries from \( D \). 

In particular, setting \( D=\mathbb {C} \), if \( u_{1},\dots ,u_{n} \) are
free from \( C \), then the matrix \( U \) is free from \( M_{N\times N}\otimes C \)
with amalgamation over the algebra of scalar diagonal matrices.
\end{lem}
The proof of the lemma can be obtained by straightforward computation of moments,
and is omitted.

\begin{lem}
\label{label:almostfree}Let \( N \) and \( s \) be fixed. Let \( A \) be
the algebra of \( N\times N \) scalar matrices, and \( B \) be the algebra
of \( dN\times dN \) scalar matrices. Denote by \( E_{\Delta } \) the \( \frac{1}{dN}\Tr  \)-preserving
conditional expectation from \( B \) onto \( \Delta _{N}\subset A=M_{N\times N}\otimes 1_{M_{d\times d}}\subset M_{N\times N}\otimes M_{d\times d}=B \).
Let \( U(d)^{\oplus N} \) denote unitaries in \( B \) which commute with \( E_{\Delta } \),
and denote by \( \mu  \) the normalized Haar measure on this compact Lie group.
Given \( \epsilon >0 \), \( \delta >0 \), \( \alpha >0 \), \( l>0 \) and
elements \( d_{1},\dots ,d_{m}\in \Delta _{N} \), there exists a \( d_{0}>0 \)
so that for all \( d>d_{0} \), given \( U_{n+1},\dots ,U_{n+r}\in U(dN) \),
there is a subset \( X\subset (U(d)^{\oplus N})^{s} \) so that \( \mu (X)>1-\alpha  \),
and so that for all \( (U_{1},\dots ,U_{s})\in X \), one has:
\[
|E_{\Delta _{N}}(d_{i_{0}}U_{j_{1}}^{g(1)}\dots d_{i_{k}}U_{j_{k}}^{g(k)})-E_{\Delta _{N}}(d_{i_{0}}u_{j_{1}}^{g(1)}\dots d_{i_{k}}u_{j_{k}}^{g(k)}|_{\epsilon }<\delta \]
for all \( k\leq l \), \( i_{1},\dots ,i_{k}\in \{1,\dots ,m\} \), \( j_{1},\dots ,j_{k}\in \{1,\dots ,n,n+1,\dots ,n+r\} \)
and \( g:\{1,\dots ,k\}\to \{\pm 1\} \). Here \( u_{1},\dots ,u_{n} \) are
Haar unitaries (all non-trivial moments are zero), which are independent from
\( \Delta _{N} \) and free with amalgamation over \( \Delta _{N} \) from each
other and from \( \{u_{n+1},\dots ,u_{n+r}\} \) (we set \( u_{j}=U_{j} \)
for \( j>n \)). In other words, \( s \)-tuples from \( X \) consists of elements
which are free among each other and from \( U_{n+1},\dots ,U_{n+r} \) with
amalgamation over \( \Delta  \) up to order \( l \) and degree \( \delta  \)
in \( |\cdot |_{\epsilon } \).
\end{lem}
The proof is a straightforward adaptation of a the proof of a similar statement
in \cite{dvv:improvedrandom}; the key observation is that because of Lemma
\ref{lemma:Ufreediag}, the desired approximate freeness holds if the \( M_{d\times d} \)-valued
entries of \( U_{1},\dots ,U_{n} \) (which are unitaries from \( U(d) \))
are \( l,\delta  \)-free from each other and also from the \( M_{d\times d} \)-valued
entries of \( U_{n+1},\dots ,U_{n+r} \) . The existence of the set \( X \)
is now guaranteed by a result in \cite{dvv:improvedrandom}.

\begin{cor}
\label{corr:gammasalmostfree}Given \( N>0 \), \( d_{1},\dots ,d_{m}\in \Delta _{N} \),
\( \epsilon ,\delta ,\alpha >0 \) and \( l>0 \), there exists a universal
constant \( d_{0} \) so that for all \( d>d_{0} \), whenever \( \Gamma _{1},\dots ,\Gamma _{n}\in U(M_{dN\times dN}) \)
are open sets, so that for each \( j \), \( \Gamma _{j} \) is invariant under
conjugation by unitaries from \( U(d)^{\oplus N} \), then there exists a subset
\( Y\subset \Gamma _{1}\times \dots \times \Gamma _{n} \), so that \( \mu (X)/\prod \mu (\Gamma _{j})>1-\alpha  \),
and such that for all \( (U_{1},\dots ,U_{n})\in Y \),
\[
|E_{\Delta _{N}}(d_{i_{0}}U_{j_{1}}^{g(1)}\dots d_{i_{k}}U_{j_{k}}^{g(k)})-E_{\Delta _{N}}(d_{i_{0}}u_{j_{1}}^{g(1)}\dots d_{i_{k}}u_{j_{k}}^{g(k)}|_{\epsilon }<\delta \]
for all \( k\leq l \), \( i_{1},\dots ,i_{k}\in \{1,\dots ,m\} \), \( j_{1},\dots ,j_{k}\in \{1,\dots ,n\} \)
and \( g:\{1,\dots ,k\}\to \{\pm 1\} \), where \( u_{j} \) has the same \( \Delta _{N} \)-valued
distribution as \( U_{j} \), and \( u_{1},\dots ,u_{n} \) are free with amalgamation
over \( \Delta _{N} \).
\end{cor}
\begin{proof}
Write
\[
\Gamma _{j}=\sqcup _{\gamma \in T_{j}}O_{\gamma },\]
where
\[
O_{\gamma }=\bigcup _{u\in U(d)^{\oplus n}}u\gamma u^{*}.\]
The Haar measure on \( \prod \Gamma _{j} \) disintegrates as \( d\mu (u)=d\mu _{O_{g}}(u)d\mu _{T}(g) \),
where \( d\mu _{O_{g}} \) is the induced Haar measure on the orbit \( O_{g} \).
For each \( g=(\gamma _{1},\dots ,\gamma _{n})\in \prod T_{j} \), let \( X \)
be the set given in Lemma \ref{corr:gammasalmostfree} for \( U_{n+1}=\gamma _{1},\dots ,U_{n+r}=\gamma _{r} \).
Let
\[
\hat{O}_{g}=\bigcup _{(u_{1},\dots ,u_{n})\in X}(u_{1}\gamma _{1}u_{1}^{*},\dots ,u_{n}\gamma _{n}u_{n}^{*}).\]
 Then \( \mu _{Og}(\hat{O}_{g})/\mu _{O_{g}}(O_{g})>1-\alpha  \). Letting \( Y=\sqcup _{g\in \prod T_{j}}\hat{O}_{g} \)
gives the statement.
\end{proof}

\section{Free Entropy \protect\( \chi (\cdots :\cdots \btw B)\protect \).}

\subsection{Sets of microstates.}

Let \( M \) be a von Neumann algebra, \( L^{\infty }[0,1]\cong B\subset M \)
a unital subalgebra, and \( u_{1},\dots ,u_{n}\in M \) be unitaries, normalizing
\( B \). Given \( \sigma =(\sigma _{1},\dots ,\sigma _{n})\in S_{N} \) and
\( d_{1},\dots ,d_{l}\in L^{\infty }[0,1] \), set
\begin{eqnarray*}
\Gamma ^{\sigma }(u_{1},\dots ,u_{n}:d_{1},\dots ,d_{m},\epsilon ,\delta ,l,d,N) & = & \\
\{(U_{1},\dots ,U_{n})\in (\sigma _{1}\cdot (U(d)^{\oplus N}),\dots ,\sigma _{n}\cdot (U(d))^{\oplus N}): &  & \\
|E_{\Delta _{N}}(d_{j_{0}}U_{i_{1}}^{g(1)}d_{j_{1}}\dots U_{i_{k}}^{g(k)}d_{j_{k}})-E_{\Delta _{N}}(d_{j_{0}}u_{i_{1}}^{g(1)}d_{j_{1}}\dots u_{i_{k}}^{g(k)}d_{j_{k}})|_{\epsilon }<\delta \} &  & 
\end{eqnarray*}
for all \( 1\leq k\leq l \), \( i_{1},\dots ,i_{k}\in \{1,\dots ,n\} \), \( j_{0},\dots ,j_{k}\in \{1,\dots ,m\} \)
and \( g:\{1,\dots ,k\}\to \{\pm 1\} \).

\begin{defn}
We shall write
\begin{eqnarray*}
\Gamma (u_{1},\dots ,u_{n}:d_{1},\dots ,d_{m},\sigma ,\epsilon ,\delta ,l,d,N)=\sigma ^{-1}\cdot \Gamma ^{\sigma }= &  & \\
\{(\sigma _{1}^{-1}U_{1},\dots ,\sigma _{n}^{-1}U_{n}):U_{1},\dots ,U_{n}\in \Gamma ^{\sigma }(u_{1},\dots ,u_{n}:d_{1},\dots ,d_{m},\epsilon ,\delta ,l,d,N) &  & 
\end{eqnarray*}
.

Define
\begin{eqnarray*}
\Gamma ^{\sigma }(u_{1},\dots ,u_{n}:u_{n+1},\dots ,u_{m}:d_{1},\dots ,d_{m},\epsilon ,\delta ,d,N)= &  & \\
\pi _{n}\Gamma ^{\sigma }(u_{1},\dots ,u_{n},u_{n+1},\dots ,u_{m}:d_{1},\dots ,d_{m},\epsilon ,\delta ,d,N), &  & 
\end{eqnarray*}
where \( \pi _{n} \) denotes the projection onto the first \( n \) components
in \( (M_{dN\times dN})^{n+m} \). 
\end{defn}

\subsection{Free entropy.}

The following definition is a straightforward adaptation of Voiculescu's definitions
of free entropy in \cite{dvv:entropy2}, \cite{dvv:entropy3}. We are dealing
with unitary elements, rather than self-adjoint ones. The appropriate modification
of Voiculescu's entropy for unitary matrices (in the absence of a subalgebra
\( B \)) was worked out in \cite{hiai-petz}.

\begin{defn}
\label{def:FDA}Assume that \( u_{1},\dots ,u_{n}\in M \) normalize \( B\cong L^{\infty }[0,1] \).
We say that \( u_{1},\dots ,u_{n},B \) have finite-dimensional approximants
(f.d.a) if for all \( D>0 \), \( \epsilon ,\delta >0 \) and \( d_{1},\dots ,d_{n}\in B \),
there are \( N>M \), so that for all \( D>0 \), there is a \( d>D \) for
which the set \( \Gamma ^{\sigma }(u_{1},\dots ,u_{n}:d_{1},\dots ,d_{m},\epsilon ,\delta ,l,d,N) \)
is non-empty for some \( \sigma \in S^{n}_{N} \).
\begin{defn}
Given a free ultrafilter \( \omega \in \beta \mathbb {N}\setminus \mathbb {N} \),
a von Neumann algebra \( M \), a unital subalgebra \( L^{\infty }[0,1]\cong B\subset M \)
and unitaries \( u_{1},\dots ,u_{n}\in M \), \( u_{n+1},\dots ,u_{q}\in M \),
normalizing \( B \), define successively:
\begin{eqnarray*}
\chi (u_{1},\dots ,u_{n}:u_{n+1},\dots ,u_{q}:d_{1},\dots ,d_{m},\epsilon ,\delta ,l,d,N)= &  & \\
\frac{1}{Nd^{2}}\sup _{\sigma \in (S_{N})^{n}}\log \mu (\Gamma (u_{1},\dots ,u_{n}:u_{n+1},\dots ,u_{q}:d_{1},\dots ,d_{m},\sigma ,\epsilon ,\delta ,l,d,N), &  & 
\end{eqnarray*}

\begin{eqnarray*}
\chi (u_{1},\dots ,u_{n}:u_{n+1},\dots ,u_{q}:d_{1},\dots ,d_{m},\epsilon ,l,N)= &  & \\
\limsup _{d\to \infty }\chi (u_{1},\dots ,u_{n}:u_{n+1},\dots ,u_{q}:d_{1},\dots ,d_{m},\epsilon ,\delta ,l,d,N), &  & 
\end{eqnarray*}

\begin{eqnarray*}
\chi ^{\omega }(u_{1},\dots ,u_{n}:u_{n+1},\dots ,u_{q}:d_{1},\dots ,d_{m},\epsilon ,l,N)= &  & \\
\lim _{d\to \omega }\chi (u_{1},\dots ,u_{n}:u_{n+1},\dots ,u_{q}:d_{1},\dots ,d_{m},\epsilon ,\delta ,l,d,N), &  & 
\end{eqnarray*}

\begin{eqnarray*}
\chi (u_{1},\dots ,u_{n}:u_{n+1},\dots ,u_{q}:d_{1},\dots ,d_{m},\epsilon ,\delta ,l)= &  & \\
\limsup _{N\to \infty }\chi (u_{1},\dots ,u_{n}:u_{n+1},\dots ,u_{q}:d_{1},\dots ,d_{m},\epsilon ,\delta ,l,N) &  & 
\end{eqnarray*}

\begin{eqnarray*}
\chi ^{\omega }(u_{1},\dots ,u_{n}:u_{n+1},\dots ,u_{q}:d_{1},\dots ,d_{m},\epsilon ,\delta ,l)= &  & \\
\lim _{N\to \omega }\chi ^{\omega }(u_{1},\dots ,u_{n}:u_{n+1},\dots ,u_{q}:d_{1},\dots ,d_{m},\epsilon ,\delta ,l,N) &  & 
\end{eqnarray*}

\begin{eqnarray*}
\chi (u_{1},\dots ,u_{n}:u_{n+1},\dots ,u_{q}:d_{1},\dots ,d_{m})= &  & \\
\inf _{l>0}\inf _{\epsilon ,\delta >0}\chi (u_{1},\dots ,u_{n}:u_{n+1},\dots ,u_{q}:d_{1},\dots ,d_{m},\epsilon ,\delta ,l) &  & 
\end{eqnarray*}

\begin{eqnarray*}
\chi ^{\omega }(u_{1},\dots ,u_{n}:u_{n+1},\dots ,u_{q}:d_{1},\dots ,d_{m})= &  & \\
\inf _{l>0}\inf _{\epsilon ,\delta >0}\chi ^{\omega }(u_{1},\dots ,u_{n}:u_{n+1},\dots ,u_{q}:d_{1},\dots ,d_{m},\epsilon ,\delta ,l) &  & 
\end{eqnarray*}

\begin{eqnarray*}
\chi (u_{1},\dots ,u_{n}:u_{n+1},\dots ,u_{q}\btw B)= &  & \\
\inf _{m>0}\inf _{d_{1},\dots ,d_{m}\in L^{\infty }[0,1]}\chi (u_{1},\dots ,u_{n}:u_{n+1},\dots ,u_{q}:d_{1},\dots ,d_{m}). &  & 
\end{eqnarray*}

\begin{eqnarray*}
\chi ^{\omega }(u_{1},\dots ,u_{n}:u_{n+1},\dots ,u_{q}\btw B)= &  & \\
\inf _{m>0}\inf _{d_{1},\dots ,d_{m}\in L^{\infty }[0,1]}\chi ^{\omega }(u_{1},\dots ,u_{n}:u_{n+1},\dots ,u_{q}:d_{1},\dots ,d_{m}), &  & 
\end{eqnarray*}
where \( \mu  \) denotes the normalized (total mass \( 1 \)) Haar measure
on \( U(d)^{\oplus N} \) (diagonal \( N\times N \) matrices with entries from
\( U(d) \)). We write simply \( \chi (u_{1},\dots ,u_{n}\btw B) \) in the
case that \( q=n \). The quantity \( \chi (u_{1},\dots ,u_{n}:u_{n+1},\dots ,u_{q}\btw B) \)
will be called \emph{free entropy of \( u_{1},\dots ,u_{n} \) in the presence
of \( u_{n+1},\dots ,u_{q} \) with respect to \( B \).}
\end{defn}
\end{defn}
In the case that \( q=\infty  \), we define \( \chi (u_{1},\dots ,u_{n}:u_{n+1},u_{n+1},\dots :d_{1},\dots ,d_{m},\epsilon ,\delta ,l) \)
to be the limit \( \liminf _{r\to \infty }\chi (u_{1},\dots ,u_{n}:u_{n+1},\dots ,u_{r}:d_{1},\dots ,d_{m},\epsilon ,\delta ,l) \),
and use that in the subsequent definitions of \( \chi  \). (Note that \( \liminf  \)
in this case is a limit). By default, \emph{we shall only deal with entropy
in the presence of a finite number of variable}, unless we explicitely state
otherwise. 

\begin{rem}
Notice that by definition \( \chi (u_{1},\dots ,u_{n}:u_{n+1},\dots ,u_{q}\btw B)\leq 0 \),
since \( \mu  \) has total mass \( 1 \).
\end{rem}
\begin{lem}
\label{lemma:replacesigma}Let \( \sigma =(\sigma _{1},\dots ,\sigma _{n}) \)
and \( \sigma '=(\sigma _{1}',\dots ,\sigma _{n}') \) be in \( S_{N}^{n} \).
Let \( d_{1},\dots ,d_{m}\in L^{\infty }[0,1] \). Assume that \( |\sigma _{i}(d_{j})-\sigma _{i}'(d_{j})|_{\alpha }<\beta  \)
for all \( 1\leq i\leq n \), \( 1\leq j\leq m \). Then
\[
\sigma '\cdot \sigma ^{-1}\cdot \Gamma ^{\sigma }(u_{1},\dots ,u_{n}:d_{1},\dots ,d_{m},\epsilon ,\delta ,l,d,N)\subset \Gamma ^{\sigma '}(u_{1},\dots ,u_{n}:d_{1},\dots ,d_{m},\epsilon +l\alpha ,\delta +l\beta ,l,d,N).\]

\begin{lem}
Let \( d_{1},\dots ,d_{m}\in B \). Assume that for each \( d\in L^{\infty }[0,1] \),
\( \epsilon ,\delta >0 \) there is polynomial \( p \) in \( d_{1},\dots ,d_{m} \)
for which \( |p(d_{1},\dots ,d_{n})-d|_{\epsilon }<\delta  \). Then \( \chi (u_{1},\dots ,u_{n}:u_{n+1},\dots ,u_{q}\btw B)=\chi (u_{1},\dots ,u_{n}:u_{n+1},\dots ,u_{q}:d_{1},\dots ,d_{m}) \).
\end{lem}
\end{lem}
\begin{proof}
One clearly has \( \chi (u_{1},\dots ,u_{n}:u_{n+1},\dots ,u_{q}\btw B)\leq \chi (u_{1},\dots ,u_{n}:u_{n+1},\dots ,u_{q}:d_{1},\dots ,d_{m}) \).
On the other hand, for \( p \) a polynomial of fixed degree \( r \),
\begin{eqnarray*}
\chi (u_{1},\dots ,u_{n}:u_{n+1},\dots ,u_{q}:d_{1},\dots ,d_{m},\epsilon ,\delta ,l)\leq  &  & \\
\chi (u_{1},\dots ,u_{n}:d_{1},\dots ,d_{m},p(d_{1},\dots ,d_{m}),\epsilon ,\delta ,[l/r]), &  & 
\end{eqnarray*}
where \( [\cdot ] \) denotes the integer part. This is because
\begin{eqnarray*}
\Gamma (u_{1},\dots ,u_{n}:u_{n+1},\dots ,u_{q}:d_{1},\dots ,d_{m},\sigma ,\epsilon ,\delta ,l,d,N)\subset  &  & \\
\Gamma (u_{1},\dots ,u_{n}:u_{n+1},\dots ,u_{q}:d_{1},\dots ,d_{m},p(d_{1},\dots ,d_{m}),\sigma ,\epsilon ,\delta ,[l/r],d,N). &  & 
\end{eqnarray*}
It follows that
\begin{eqnarray*}
\chi (u_{1},\dots ,u_{n}:u_{n+1},\dots ,u_{q}:d_{1},\dots ,d_{m},\epsilon ,\delta ,l)\leq  &  & \\
\chi (u_{1},\dots ,u_{n}:u_{n+1},\dots ,u_{q}:d_{1},\dots ,d_{m},d,\epsilon +2\epsilon '[l/r],\delta +2\delta '[l/r],[l/r]). &  & 
\end{eqnarray*}
 and \( |p(d_{1},\dots ,d_{n})-d|_{\epsilon '}<\delta ' \). It follows after
taking limits that 
\[
\inf _{l>0}\inf _{\epsilon ,\delta >0}\chi (u_{1},\dots ,u_{n}:d_{1},\dots ,d_{m},\epsilon ,\delta ,l)\leq \inf _{\epsilon ,\delta >0}\chi (u_{1},\dots ,u_{n}:d_{1},\dots ,d_{m},d,\epsilon ,\delta ,l)\]
which in turn implies that
\[
\chi (u_{1},\dots ,u_{n}:u_{n+1},\dots ,u_{q}:d_{1},\dots ,d_{m})\leq \chi (u_{1},\dots ,u_{n}:u_{n+1},\dots ,u_{q}:d_{1},\dots ,d_{m},d).\]

Hence whenever \( d_{1}',\dots ,d_{m'}'\in B \), we get
\begin{eqnarray*}
\chi (u_{1},\dots ,u_{n}:u_{n+1},\dots ,u_{q}:d_{1},\dots ,d_{m})\leq  &  & \\
\chi (u_{1},\dots ,u_{n}:u_{n+1},\dots ,u_{q}:d_{1},\dots ,d_{m},d_{1}',\dots ,d_{m'}')\leq  &  & \\
\chi (u_{1},\dots ,u_{n}:u_{n+1},\dots ,u_{q}:d_{1}',\dots ,d_{m'}'), &  & 
\end{eqnarray*}
which in turn gives \( \chi (u_{1},\dots ,u_{n}:u_{n+1},\dots ,u_{q}\btw B)\geq \chi (u_{1},\dots ,u_{n}:u_{n+1},\dots ,u_{q}:d_{1},\dots ,d_{m}) \).
\end{proof}
\begin{prop}
Let \( u_{1},\dots ,u_{n},v_{1},\dots ,v_{q}\in M \). Then
\begin{eqnarray*}
\chi (u_{1},\dots ,u_{n}:v_{1},\dots ,v_{q}\btw B)\leq  &  & \\
\chi (u_{1},\dots ,u_{r}:v_{1},\dots ,v_{q}\btw B)+\chi (u_{r+1},\dots ,u_{n}:v_{1},\dots ,v_{q}\btw B) &  & 
\end{eqnarray*}
and similarly for \( \chi ^{\omega } \).
\end{prop}
\begin{proof}
This follows from the obvious inclusion
\begin{eqnarray*}
\Gamma (u_{1},\dots ,u_{n},v_{1},\dots ,v_{q}:d_{1},\dots ,d_{m},\sigma ,\epsilon ,\delta ,l,d,N)\subset  &  & \\
\Gamma (u_{1},\dots ,u_{r},v_{1},\dots ,v_{q}:d_{1},\dots ,d_{m},\sigma ',\epsilon ,\delta ,l,d,N) &  & \\
\times \Gamma (u_{1},\dots ,u_{r},v_{1},\dots ,v_{q}:d_{1},\dots ,d_{m},\sigma '',\epsilon ,\delta ,l,d,N), &  & 
\end{eqnarray*}
where \( S_{N}^{n+2q}\ni \sigma =(\sigma ',\sigma '')\in S_{N}^{r+q}\times S_{N}^{n-r+q} \).
\end{proof}
\begin{prop}
\label{prop:changeofvars}Let \( u_{1},\dots ,u_{n},v_{1},\dots ,v_{q}\in M \).
Let \( w_{1},\dots ,w_{r}\in W^{*}(B,v_{1},\dots ,v_{q},u_{1},\dots ,u_{n}) \).
Then
\[
\chi (u_{1},\dots ,u_{n}:v_{1},\dots ,v_{q}\btw B)=\chi (u_{1},\dots ,u_{n}:v_{1},\dots ,v_{q},w_{1},\dots ,w_{r}\btw B)\]
and similarly for \( \chi ^{\omega } \) instead of \( \chi  \).
\end{prop}
The proof of this Proposition is essentially identical to the proof of Proposition
a similar Proposition in \cite{dvv:entropy3}, and is therefore omitted.

\begin{prop}
\( \chi (u_{1},\dots ,u_{n}:v_{1},\dots ,v_{r}\btw B)\leq \chi (u_{1},\dots ,u_{n}:v_{1},\dots ,v_{q}\btw B) \)
if \( r\leq q \).
\end{prop}
\begin{proof}
One has
\begin{eqnarray*}
\Gamma (u_{1},\dots ,u_{n},v_{1},\dots ,v_{q}:d_{1},\dots ,d_{m},\sigma ,\epsilon ,\delta ,l,d,N)\subset  &  & \\
\Gamma (u_{1},\dots ,u_{n},v_{1},\dots ,v_{r}:d_{1},\dots ,d_{m},\sigma ,\epsilon ,\delta ,l,d,N); &  & 
\end{eqnarray*}
the inequality now follows after taking limits.
\end{proof}

\section{Properties of \protect\( \chi (\cdots :\cdots \btw B)\protect \).}

\begin{prop}
\label{Prop:independence}Let \( u_{1},\dots ,u_{q}\in M \) be such that \( [u_{j},B]=\{0\} \),
and \( W^{*}(u_{1},\dots ,u_{q}) \) is independent from \( B \). Then 
\[
\chi (u_{1},\dots ,u_{n}:u_{n+1},\dots ,u_{q}\btw B)=\chi (u_{1},\dots ,u_{n}:u_{n+1},\dots ,u_{q}),\]
where the last quantity is the unitary analog of Voiculescu's entropy in the
presence (see \cite{dvv:entropy6}, \cite{hiai-petz}). The same statement holds
true for \( \chi ^{\omega } \) instead of \( \chi  \).
\end{prop}
\begin{proof}
We shall first prove that \( \chi (u_{1},\dots ,u_{n}:u_{n+1},\dots ,u_{q}\btw B)\geq \chi (u_{1},\dots ,u_{n}:u_{n+1},\dots ,u_{q}) \).
Fix \( d_{1},\dots ,d_{m}\in B \). Let \( \sigma =(\id ,\dots ,\id )\in S_{N}^{n} \),
and consider the set
\[
X=\Gamma (u_{1},\dots ,u_{n},u_{n+1},\dots ,u_{q};l,d,\delta )^{\oplus N}\subset \sigma \cdot (U(d)^{\oplus N})^{q}.\]
 We claim that \( X\subset \Gamma ^{\sigma }(u_{1},\dots ,u_{q}:d_{1},\dots ,d_{m},\epsilon ,\delta ,l,d,N) \).
To show this, it is sufficient to verify (by enlarging the set \( d_{1},\dots ,d_{n} \)
to contain all words in \( d_{1},\dots ,d_{n} \) of length at most \( l \)
and also the unit of \( B \)) that for \( \Vert d_{j}\Vert _{\infty }\leq 1 \),
\( 1\leq j\leq m \), and for \( k\leq l \) and \( i_{1},\dots ,i_{k}\in \{1,\dots ,q\} \),
\( g:\{1,\dots ,k\}\to \{\pm 1\} \),
\[
\left| d_{j}\left[ E_{\Delta _{N}}(u_{i_{1}}^{g(1)}\dots u_{i_{k}}^{g(k)})-E_{\Delta _{N}}(U_{i_{1}}^{g(1)}\dots U_{i_{k}}^{g(k)})\right] \right| _{\epsilon }<\delta ,\]
or, equivalently, using independence of \( W^{*}(u_{1},\dots ,u_{n}) \) and
\( B \),
\[
\left| d_{j}\left[ \tau (u_{i_{1}}^{g(1)}\dots u_{i_{k}}^{g(k)})-E_{\Delta _{N}}(U_{i_{1}}^{g(1)}\dots U_{i_{k}}^{g(k)})\right] \right| _{\epsilon }<\delta ,\]
for all \( (U_{1},\dots ,U_{q})\in X \). Writing \( U_{j}=w^{(1)}_{j}\oplus \dots \oplus w_{j}^{(N)} \),
we see that the equation above is satisfied if \( \tau (u_{i_{1}}^{g(1)}\dots u_{i_{k}}^{g(k)})-\frac{1}{d}\Tr ((w_{i_{1}}^{(j)})^{g(1)}\dots (w_{i_{k}}^{(j)})^{g(k)})<\delta  \)
for all \( k\leq l \), \( i_{1},\dots ,i_{k}\in \{1,\dots ,q\} \) and \( g:\{1,\dots ,k\}\to \{\pm 1\} \).
But this is precisely the condition that \( (U_{1},\dots ,U_{n})\in X \).

It follows that
\[
\frac{1}{Nd^{2}}\log \mu (\pi _{n}(X))\leq \chi (u_{1},\dots ,u_{n}:u_{n+1},\dots ,u_{q}:d_{1},\dots ,d_{m},\epsilon ,\delta ,l,d,N);\]
since \( \pi _{n}X=\Gamma (u_{1},\dots ,u_{n}:u_{n+1},\dots ,u_{q};l,d,\epsilon )^{\oplus N} \),
we get that
\begin{eqnarray*}
\frac{1}{Nd^{2}}\mu (X)=\frac{1}{d^{2}}\frac{1}{N}\log \mu (\Gamma (u_{1},\dots ,u_{n}:u_{n+1},\dots ,u_{q};l,d,\epsilon )^{\oplus N})= &  & \\
\chi (u_{1},\dots ,u_{n}:u_{n+1},\dots ,u_{q};l,d,\epsilon )\leq  &  & \\
\chi (u_{1},\dots ,u_{n}::u_{n+1},\dots ,u_{q}:d_{1},\dots ,d_{m},\epsilon ,\delta ,l,d,N), &  & 
\end{eqnarray*}
which implies \( \chi (u_{1},\dots ,u_{n}:u_{n+1},\dots ,u_{q}\btw B)\geq \chi (u_{1},\dots ,u_{n}:u_{n+1},\dots ,u_{q}) \).

To prove the opposite inequality, let now \( \sigma ' \) be such that 
\begin{eqnarray*}
\mu (\Gamma ^{\sigma '}(u_{1},\dots ,u_{n}:u_{n+1},\dots ,u_{q}:d_{1},\dots ,d_{m},\epsilon ,\delta ,l,d,N))= &  & \\
\sup _{\sigma ''\in S_{N}^{n}}\mu (\Gamma ^{\sigma ''}(u_{1},\dots ,u_{n}:u_{n+1},\dots ,u_{q}:d_{1},\dots ,d_{m},\epsilon ,\delta ,l,d,N)). &  & 
\end{eqnarray*}
Then
\[
|\sigma _{i}'(d_{j})-\sigma _{i}(d)|_{\epsilon }<\delta \]
for all \( i,j \). Hence by Lemma \ref{lemma:replacesigma},
\begin{eqnarray}
\sigma \cdot (\sigma ')^{-1}\Gamma ^{\sigma '}(u_{1},\dots ,u_{n}:u_{n+1},\dots ,u_{q}:d_{1},\dots ,d_{m},\epsilon ,\delta ,l,d,N) & \subset  & \label{eqn:sigmatoidentity} \\
\Gamma ^{\sigma }(u_{1},\dots ,u_{n}:u_{n+1},\dots ,u_{q}:d_{1},\dots ,d_{m},\epsilon (1+l),\delta (1+l),l,d,N) &  & 
\end{eqnarray}
Since the unit of \( B \) occurs among \( d_{1},\dots ,d_{n} \), this implies
that any
\[
(U_{1},\dots ,U_{q})\in \Gamma ^{\sigma }(u_{1},\dots ,u_{n},u_{n+1},\dots ,u_{q}:d_{1},\dots ,d_{m},\epsilon (1+l),\delta (1+l),l,d,N)\]
satisfy
\[
|E_{\Delta _{N}}(U_{j_{1}}^{g(1)}\dots U_{j_{k}}^{g(k)})-\tau (u_{j_{1}}^{g(1)}\dots u_{j_{k}}^{g(k)})|_{\epsilon (1+l)}<\delta (1+l),\]
for all \( k\leq l \), \( j_{1},\dots ,j_{k}\in \{1,\dots ,n\} \), \( g:\{1,\dots ,k\}\to \{\pm 1\} \).
Let
\[
U_{j}=w^{(1)}_{j}\oplus \dots \oplus w_{j}^{(N)}\]
with \( w_{j}^{(k)}\in U(d) \). Notice that 
\[
E_{\Delta _{N}}(U_{j_{1}}^{g(1)}\dots U_{j_{k}}^{g(k)})\]
is a diagonal matrix, whose \( r \)-th diagonal entry is
\[
\frac{1}{d}\Tr ((w_{j_{1}}^{(r)})^{g(1)}\dots (w_{j_{k}}^{(r)})^{g(k)}).\]
If \( N \) is so large that \( M/N>\epsilon (1+l)n^{l} \) for some \( M<N \),
it follows that for each \( (U_{1},\dots ,U_{n}) \) there exists a subset \( S \)
of \( \{1,\dots ,N\} \) with \( |S|>N-M \), and so that for all \( r\in S \),
\[
(w_{1}^{(r)},\dots ,w_{n}^{(r)})\in \Gamma (u_{1},\dots ,u_{q};l,d,\delta (1+l)).\]
Hence
\[
\Gamma ^{\sigma }(u_{1},\dots ,u_{n}:u_{n+1},\dots ,u_{q}:d_{1},\dots ,d_{m},\epsilon (1+l),\delta (1+l),l,d,N)\subset \bigcup _{\begin{array}{c}
S\subset \{1,\dots ,N\}\\
|S|>N-M
\end{array}}\bigoplus _{p=1}^{N}X(p,S)\]
where \( X(p,S)=\Gamma (u_{1},\dots ,u_{n}:u_{n+1},\dots ,u_{q};l,d,\delta (1+l)) \)
if \( p\in S \) and \( X(p,S)=U(d) \) if \( p\notin S \). It follows that
\begin{eqnarray*}
\frac{1}{Nd^{2}}\log \mu (\Gamma ^{\sigma }(u_{1},\dots ,u_{n}:d_{1},\dots ,d_{m},\epsilon (1+l),\delta (1+l),l,d,N)) & \leq  & \\
\frac{N-M}{N}\chi (u_{1},\dots ,u_{n};l,d,\delta (1+l)+ &  & \\
\frac{M}{Nd^{2}}\log (1)+\frac{1}{Nd^{2}}\log \left( \begin{array}{c}
M\\
N
\end{array}\right)  &  & 
\end{eqnarray*}
Taking the limit \( d\to \infty  \) and using (\ref{eqn:sigmatoidentity})
gives
\begin{eqnarray*}
\chi (u_{1},\dots ,u_{n}:u_{n+1},\dots ,u_{q}:d_{1},\dots ,d_{m},\epsilon ,\delta ,l,N)\leq  &  & \\
\frac{N-M}{N}\chi (u_{1},\dots ,u_{n}:u_{n+1},\dots ,u_{q};l,d,\delta (l+1))+ &  & \\
\lim _{d\to \infty }\frac{1}{Nd^{2}}\log \left( \begin{array}{c}
M\\
N
\end{array}\right) . &  & 
\end{eqnarray*}
Since \( M \) is chosen so that \( M/N>\epsilon (1+l)n^{l} \), taking the
limit as \( N\to \infty  \) and infimum over \( \epsilon ,\delta  \) and \( l \)
gives the desired inequality. 

The proof for \( \chi ^{\omega } \) is identical.
\end{proof}
\begin{prop}
\label{prop:estimateonchi}Let \( u\in M \) be a unitary, so that \( [u,B]=0 \).
Assume that \( \Vert E_{B}(|u-1|^{2})\Vert ^{1/2}_{\infty }<\delta  \). Then
\[
\chi (u\btw B)\leq \log \delta +C,\]
for some universal constant \( C \).
\end{prop}
\begin{proof}
Let \( d_{1},\dots ,d_{m}\in B \), \( \Vert d_{j}\Vert _{\infty }\leq 1 \)
be given. Let \( \epsilon >0 \). Then by Lemma \ref{lemma:replacesigma} we
have that \( \Gamma (u:d_{1},\dots ,d_{m},\sigma ,\epsilon /2,\delta ^{2}/4,2,d,N)\subset \Gamma (u:d_{1},\dots ,d_{m},\id ,\epsilon ,\delta ^{2},2,d,N)=\Gamma  \).
Let \( U=U_{1}\oplus \dots \oplus U_{N}\in \Gamma  \). Then we have in particular
that
\[
|E_{\Delta _{N}}((U-I)(U-I)^{*})-E_{\Delta _{N}}((u-1)(u-1)^{*})|_{\epsilon }<4\delta ^{2}.\]
Note that \( \Vert E_{\Delta _{N}}((u-1)(u-1)^{*})\Vert _{\infty }\leq \Vert E_{B}(|u-1|^{2})\Vert ^{2}<\delta ^{2} \).
Let \( M=[\epsilon N] \), where \( [\cdot ] \) denotes the integer part of
a number. Then for at least \( N-M \) numbers \( j \) in the set \( \{1,\dots ,N\} \),
we have \( \Vert U_{j}-1\Vert ^{2}_{2}\leq 5\delta ^{2}<(3\delta )^{2} \).
It follows that \( \Gamma  \) is contained in the set
\[
\Gamma \subset S=\bigsqcup _{J\subset \{1,\dots ,N\},|J|>\epsilon N}S(1,J)\oplus \dots \oplus S(N,J),\]
where \( S(k,J)=U(d) \) when \( k\notin J \) and \( S(k,J) \) is the ball
\( S(k,J)=B(U(d),3\delta )=\{U\in U(d):\Vert U-I\Vert _{2}\leq 3\delta \} \)
for \( k\in J \) (the \( \Vert \cdot \Vert _{2} \) norm is with respect to
the normalized trace \( \frac{1}{d}\Tr  \) on \( U(d) \)). It follows that
\[
\frac{1}{Nd^{2}}\log \mu (\Gamma )\leq \frac{1}{Nd^{2}}\log \left( \begin{array}{c}
N\\
N-[\epsilon N]
\end{array}\right) +\frac{N-[\epsilon N]}{Nd^{2}}\log \mu (B(U(d),3\delta )).\]
The limit as \( d\to \infty  \) of \( \frac{1}{Nd^{2}}\log \left( \begin{array}{c}
N\\
N-[\epsilon N]
\end{array}\right)  \) is zero. Hence we get
\[
\chi (u_{1},\dots ,u_{n}:d_{1},\dots ,d_{m},\epsilon /2,\delta ^{2}/4,2,N)\leq \lim _{d}\frac{N-[\epsilon N]}{N}\frac{1}{d^{2}}\log \mu (B(U(d),3\delta )).\]
 As \( d\to \infty  \) and \( N\to \infty  \), we get as estimate \( (1-\epsilon )\log \delta +C \)
for some universal constant \( C \). The desired estimate now follows from
the definition of \( \chi  \).
\end{proof}
\begin{prop}
\label{prop:additivityundercompression}Let \( u_{1},\dots ,u_{q}\in M \),
and assume that \( p_{1},\dots ,p_{r}\in B\cong L^{\infty }[0,1] \) are projections,
\( \sum p_{i}=1 \), so that \( [u_{i},p_{j}]=0 \) for all \( 1\leq i\leq n \)
and \( 1\leq j\leq r \), with possibly \( r=\infty  \). Then \( p_{j}u_{i}p_{j} \)
is a unitary in the algebra \( p_{j}Mp_{j} \), and \( L^{\infty }[0,1]\cong p_{j}Bp_{j}\subset p_{j}Mp_{j} \). 

We have
\begin{eqnarray*}
\chi ^{\omega }(u_{1},\dots ,u_{n}:u_{n+1},\dots ,u_{q}\btw B)= &  & \\
\sum _{j=1}^{r}\tau (p_{j})\chi ^{\omega \cdot \tau (p_{j})}(p_{j}u_{1}p_{j},\dots ,p_{j}u_{n}p_{j}:p_{j}u_{n+1}p_{j},\dots ,p_{j}u_{q}p_{j}\btw p_{j}Bp_{j}). &  & 
\end{eqnarray*}
and
\begin{eqnarray*}
\chi (u_{1},\dots ,u_{n}:u_{n+1},\dots ,u_{q}\btw B)= &  & \\
\sum _{j=1}^{r}\tau (p_{j})\chi (p_{j}u_{1}p_{j},\dots ,p_{j}u_{n}p_{j}:p_{j}u_{n+1}p_{j},\dots ,p_{j}u_{q}p_{j}\btw p_{j}Bp_{j}). &  & 
\end{eqnarray*}
Here \( \omega \cdot t \) for \( t\in \mathbb {R}_{+} \)denotes the ultrafilter
determined by
\[
\lim _{n\to \omega \cdot t}f(n)=\lim _{t\to \omega }f([nt]),\]
where \( [x] \) denotes the integer part of \( x \), and \( f \) is a bounded
real function on \( \mathbb {N} \).
\end{prop}
\begin{proof}
We may identify \( B \) with \( L^{\infty }[0,1] \) in such a way that the
projections \( p_{j} \) correspond to characteristic functions of the intervals
\( [x_{j},x_{j+1}] \) for some points \( 0=x_{1}\leq x_{2}\leq \dots \leq x_{r}\leq x_{r+1}=1 \).
Fix \( d_{1},\dots ,d_{n}\in B \); we can choose \( d_{1},\dots ,d_{m} \)
in such a way that \( p_{j}d_{i}=0 \) for all \( 1\leq i\leq m \) and all
\( j>j_{0} \). Choose integers \( N_{1},\dots ,N_{r} \) so that \( N_{j} \)
are zero starting from some \( j_{0} \), and and write \( d_{s}^{(j)}=p_{j}d_{s}p_{j} \),
\( N=\sum _{j=1}^{r}N_{j} \) (note that \( N_{j} \) are zero for sufficiently
large \( j \)). Then choosing \( \sigma ^{(j)}\in S_{N_{j}}^{n} \) and letting
\( \sigma =\bigoplus \sigma ^{(j)}\in S_{N}^{n} \), we have that
\begin{eqnarray*}
\bigoplus _{j=1}^{j_{0}}\Gamma (p_{j}u_{1}p_{j},\dots ,p_{j}u_{n}p_{j}:p_{j}u_{n+1}p_{j},\dots ,p_{j}u_{q}p_{j}:d_{1}^{(j)},\dots ,d_{m}^{(j)},\sigma ^{(j)},\epsilon ,\delta ,l,d,N_{j})\subset  &  & \\
\Gamma (u_{1},\dots ,u_{n}:u_{n+1},\dots ,u_{q}:d_{1},\dots ,d_{m},\sigma ,\epsilon ',\delta ,l,N), &  & 
\end{eqnarray*}
provided that \( \epsilon '\leq \sum \epsilon _{j}\frac{N_{j}}{N}+\alpha _{j} \),
where \( \alpha _{j} \) is the Lebesgue measure of the symmetric difference
of \( [x_{j},x_{j+1}] \) and \( [\frac{\sum _{i<j}N_{i}}{N},\frac{\sum _{i\leq j}N_{i}}{N}] \).
Hence for \( N \) sufficiently large, we can choose \( N_{j}=[(x_{j+1}-x_{j})N] \)
for \( 1\leq j<r \), \( j_{0} \) to be the first \( j \) for which \( N_{j} \)
is zero and \( N_{r}=N-\sum _{1\leq j<j_{0}}N_{j} \), and have that:
\begin{eqnarray*}
\bigoplus _{j=1}^{j_{0}}\Gamma (p_{j}u_{1}p_{j},\dots ,p_{j}u_{n}p_{j}:p_{j}u_{n+1}p_{j},\dots ,p_{j}u_{q}p_{j}:d_{1}^{(j)},\dots ,d_{m}^{(j)},\sigma ^{(j)},\epsilon ,\delta ,l,d,N_{j})\subset  &  & \\
\Gamma (u_{1},\dots ,u_{n}:u_{n+1},\dots ,u_{q}:d_{1},\dots ,d_{m},\sigma ,2\epsilon ,\delta ,l,N). &  & 
\end{eqnarray*}
This implies that
\begin{eqnarray*}
\sum _{j=1}^{r}\frac{N_{j}}{N}\chi (p_{j}u_{1}p_{j},\dots ,p_{j}u_{n}p_{j}:p_{j}u_{n+1}p_{j},\dots ,p_{j}u_{q}p_{j}:d_{1}^{(j)},\dots ,d_{n}^{(j)},\epsilon ,\delta ,l,d,N_{j})\leq  &  & \\
\chi (u_{1},\dots ,u_{n}:u_{n+1},\dots ,u_{q}:d_{1},\dots ,d_{m},2\epsilon ,\delta ,l,d,N). &  & 
\end{eqnarray*}
 Taking the limit \( N\to \omega  \) and noticing that in this case each \( N_{j}\to \omega \cdot \tau (p_{j}) \),
since \( \tau (p_{j})=x_{j+1}-x_{j} \), and \( N_{j}/N\to \tau (p_{j}) \),
gives that
\begin{eqnarray*}
\chi ^{\omega }(u_{1},\dots ,u_{n}:u_{n+1},\dots ,u_{q}\btw B)\geq  &  & \\
\sum _{j=1}^{r}\tau (p_{j})\chi ^{\omega \cdot \tau (p_{j})}(p_{j}u_{1}p_{j},\dots ,p_{j}u_{n}p:p_{j}u_{n+1}p_{j},\dots ,p_{j}u_{q}p_{j}\btw p_{j}Bp_{j}). &  & 
\end{eqnarray*}

Note that we have the same inequality for \( \chi ^{\omega } \) and \( \chi ^{\omega \cdot \tau (p_{j})} \)
replaced by \( \chi  \).

For the opposite inequality, we may assume that \( N=\sum _{j=1}^{j_{0}}N_{j}+k \),
with \( |N_{j}-[(x_{j+1}-x_{j})N]|\leq 1 \), \( j_{0}\leq r \) and \( k/N<\epsilon /2 \).
Moreover, assume that for some \( \sigma \in S_{N} \),
\begin{eqnarray*}
\chi (u_{1},\dots ,u_{n}:u_{n+1},\dots ,u_{q}:d_{1},\dots ,d_{n},\epsilon ,\delta ,l,d,N)= &  & \\
\frac{1}{d^{2}N}\log \mu (\Gamma (u_{1},\dots ,u_{n}:u_{n+1},\dots ,u_{q}:d_{1},\dots ,d_{n},\sigma ,\epsilon ,\delta ,l,d,N). &  & 
\end{eqnarray*}
Since \( [u_{i},p_{j}]=0 \), it follows that, given \( \epsilon '>0 \) we
can find \( \sigma ^{(1)}\in S_{N_{1}}^{n},\dots ,\sigma ^{(r)}\in S_{N_{j_{0}}}^{n} \)
and \( \epsilon >0 \) (independent of \( N \) and \( d \)), for which, after
letting \( \sigma '=\bigoplus \sigma ^{(j)}\oplus \id _{k}\in S_{\sum N_{j}+k}^{n}=S_{N}^{n} \),
one has
\begin{eqnarray*}
\Gamma (u_{1},\dots ,u_{n}:u_{n+1},\dots ,u_{q}:d_{1},\dots ,d_{n},\sigma ',\epsilon ,\delta ,l,d,N)\supset  &  & \\
\Gamma (u_{1},\dots ,u_{n}:u_{n+1},\dots ,u_{q}:d_{1},\dots ,d_{n},\sigma ,\epsilon ',\delta ,l,d,N). &  & 
\end{eqnarray*}
Let now
\[
(U_{1},\dots ,U_{n},U_{n+1},\dots ,U_{q})\in \Gamma (u_{1},\dots ,u_{n}:u_{n+1},\dots ,u_{q}:d_{1},\dots ,d_{n},\sigma ',\epsilon ,\delta ,l,d,N).\]
Let \( M=[N\epsilon /2]+1 \). Denote by \( P_{j} \) the diagonal matrix having
all entries zero, except that the \( k,k \)-th entries for \( N_{j}\leq k<N_{j+1} \)
are equal to \( 1 \). Then for a subset \( S\subset \{1,\dots ,N\} \) of size
at most \( M \), we have that
\[
(P_{j}U_{1}P_{J},\dots ,P_{j}U_{q}P_{j})\in \Gamma (p_{j}u_{1}p_{j},\dots ,p_{j}u_{q}p_{j}:p_{j}d_{1}p_{j},\dots ,p_{j}d_{n}p_{j},\sigma ^{(j)},\frac{N}{N_{j}}\epsilon ,\delta ,l,d,N_{j}).\]
Therefore, one has
\begin{eqnarray*}
\sum _{j=1}^{r}\frac{N_{j}}{N}\chi (p_{j}u_{1}p_{j},\dots ,p_{j}u_{n}p_{j}:p_{j}u_{n+1}p_{j},\dots ,p_{j}u_{q}p_{j}:d_{1}^{(j)},\dots ,d_{n}^{(j)},2\tau (p_{j})^{-1}\epsilon ,\delta ,l,d,N_{j})\geq  &  & \\
\chi (u_{1},\dots ,u_{n}:u_{n+1},\dots ,u_{q}:d_{1},\dots ,d_{m},\epsilon ,\delta ,l,d,N)-\frac{1}{Nd^{2}}\log \left( \begin{array}{c}
N\\
M
\end{array}\right) . &  & 
\end{eqnarray*}
Taking the limits \( N\to \omega  \) (so that \( N_{j}\to \tau (p_{j})\omega  \))
gives finally
\begin{eqnarray*}
\chi ^{\omega }(u_{1},\dots ,u_{n}:u_{n+1},\dots ,u_{q}\btw B)\leq  &  & \\
\sum _{j=1}^{r}\tau (p_{j})\chi ^{\omega \cdot \tau (p_{j})}(p_{j}u_{1}p_{j},\dots ,p_{j}u_{n}p_{j}:p_{j}u_{n+1}p_{j},\dots ,p_{j}u_{q}p_{j}\btw p_{j}Bp_{j}). &  & 
\end{eqnarray*}
Note that the same argument gives the same inequality for \( \chi  \) instead
of \( \chi ^{\omega } \).
\end{proof}
\begin{prop}
\label{prop:inpresencefree}Assume that \( v_{1},\dots ,v_{r} \) are free with
amalgamation over \( B \) from \( u_{1},\dots ,u_{q} \). Assume that \( B,v_{1},\dots ,v_{r} \)
has f.d.a (see Definition \ref{def:FDA}). Then
\[
\chi (u_{1},\dots ,u_{n}:u_{n+1},\dots ,u_{q},v_{1},\dots ,v_{r})=\chi (u_{1},\dots ,u_{n}:u_{n+1},\dots ,u_{q})\]
and similarly for \( \chi ^{\omega } \). The same conclusion holds for \( r=\infty  \).
\end{prop}
\begin{proof}
Fix
\[
(v_{1},\dots ,v_{r})\in \Gamma ^{\sigma }(v_{1},\dots ,v_{r}:d_{1},\dots ,d_{n},\epsilon ,\delta ,l,d,N).\]
By \ref{corr:gammasalmostfree}, for all \( \alpha >0 \), there exist a subset
\( W \)
\[
W\subset \Gamma ^{\sigma }(u_{1},\dots ,u_{n}:u_{n+1},\dots ,u_{q}:d_{1},\dots ,d_{n},\epsilon ,\delta ,l,d,N),\]
so that
\[
\frac{1}{k^{2}N}\mu (W)/\mu (\Gamma ^{\sigma }(u_{1},\dots ,u_{n}:u_{n+1},\dots ,u_{q}:d_{1},\dots ,d_{n},\epsilon ,\delta ,l,d,N))>1-\alpha ,\]
and such that if \( (w_{1},\dots ,w_{n})\in W \), then there exist \( (w_{n+1},\dots ,w_{q}) \)
so that
\[
(w_{1},\dots ,w_{n},w_{n+1},\dots ,w_{q})\in \Gamma ^{\sigma }(u_{1},\dots ,u_{n},u_{n+1},\dots ,u_{q}:d_{1},\dots ,d_{n},\epsilon ,\delta ,l,d,N).\]
and \( (w_{1},\dots ,w_{q}) \) is free up to order \( l \) and degree \( \delta  \)
in \( |\cdot |_{\epsilon } \) from \( (v_{1},\dots ,v_{r}) \) with amalgamation
over \( \Delta _{N} \). This implies that
\[
W\subset \Gamma ^{\sigma }(u_{1},\dots ,u_{n}:u_{n+1},\dots ,u_{q},v_{1},\dots ,v_{r}:d_{1},\dots ,d_{n},\epsilon ,\delta ,l,d,N).\]
Passing to the limit gives
\[
\chi (u_{1},\dots ,u_{n}:u_{n+1},\dots ,u_{q},v_{1},\dots ,v_{r})\leq \chi (u_{1},\dots ,u_{n}:u_{n+1},\dots ,u_{q}).\]
The reverse inequality is obvious.
\end{proof}
\begin{prop}
\label{prop:freeadditivity}Let \( u_{1},\dots ,u_{n},v_{1},\dots ,v_{n}\in M \),
and let \( 1<s<n \). Assume that the sets \( (u_{1},v_{1},\dots ,u_{s},v_{s}),\dots ,(u_{s+1},v_{s+1},\dots u_{n},v_{n}) \)
are \( * \)-free with amalgamation over \( B \). Then
\[
\chi ^{\omega }(u_{1},\dots ,u_{n}:v_{1},\dots ,v_{n}\btw B)=\chi ^{\omega }(u_{1},\dots ,u_{s}:v_{1},\dots ,v_{s}\btw B)+\chi ^{\omega }(u_{s+1},\dots ,u_{n}:v_{s+1},\dots ,v_{n}\btw B).\]

\end{prop}
\begin{proof}
Note first that because of the freeness assumptions,
\[
\chi ^{\omega }(u_{1},\dots ,u_{s}:v_{1},\dots ,v_{s}\btw B)=\chi ^{\omega }(u_{1},\dots ,u_{s}:v_{1},\dots ,v_{n}\btw B)\]
and
\[
\chi ^{\omega }(u_{s+1},\dots ,u_{n}:v_{1},\dots ,v_{n}\btw B).\]
The inequality
\[
\chi (u_{1},\dots ,u_{n}:v_{1},\dots ,v_{n}\btw B)\leq \chi ^{\omega }(u_{1},\dots ,u_{s}:v_{1},\dots ,v_{s}\btw B)+\chi ^{\omega }(u_{s+1},\dots ,u_{n}:v_{s+1},\dots ,v_{n}\btw B)\]
is then clear.

Fix \( N,d,l,\epsilon ,\delta ,d_{1},\dots ,d_{m} \). Choose \( \sigma _{1},\dots ,\sigma _{n}\in S_{N} \)
so that for each \( j \),
\[
\mu (\Gamma ^{\sigma _{j}}(u_{j}:v_{j}:d_{1},\dots ,d_{m},\epsilon ,\delta ,l,d,N))=\sup _{\sigma '\in S_{N}}\mu (\Gamma ^{\sigma '}(u_{j}:v_{j}:d_{1},\dots ,d_{m},\epsilon ,\delta ,l,d,N)).\]
 By \ref{corr:gammasalmostfree}, for all \( \alpha >0 \), there exist a subset
\( W \)
\begin{eqnarray*}
W\subset \Gamma ^{\sigma _{1}}(u_{1},\dots ,u_{s},v_{1},\dots ,v_{n}:d_{1},\dots ,d_{m},\epsilon ,\delta ,l,d,N)\times  &  & \\
\Gamma ^{\sigma _{2}}(u_{s+1},\dots ,u_{n},v_{1},\dots ,v_{n}:d_{1},\dots ,d_{m},\epsilon ,\delta ,l,d,N)=\Gamma  &  & 
\end{eqnarray*}
so that
\[
\frac{1}{d^{2}N}\mu (W)/\mu (\Gamma )>1-\alpha \]
and such that if \( ((W_{1},V_{1}),\dots ,(W_{n},V_{n}))\in W \), then \( (W_{1},V_{1},\dots ,W_{s},V_{s}) \)
and \( (W_{s+1},V_{s+1},\dots ,W_{n},V_{n}) \) are free up to order \( l \)
and degree \( \delta  \) in \( |\cdot |_{\epsilon } \) with amalgamation over
\( \Delta _{N} \). It follows that
\[
W\subset \mu (\Gamma ^{\sigma _{1}\oplus \sigma _{2}}(u_{1},\dots ,u_{n}:v_{1},\dots ,v_{n}:d_{1},\dots ,d_{m},\epsilon ,\delta ,l,d,N))\]
which implies the proposition after taking limits.
\end{proof}
We don't know if the preceding proposition holds for \( \chi  \) instead of
\( \chi ^{\omega } \), because there is no guarantee that the \( \limsup _{d} \)
and \( \limsup _{N} \) in the definitions of \( \chi (u_{1},\dots ,u_{s}\btw B) \)
and \( \chi (u_{s+1},\dots ,u_{n}\btw B) \) are attained on the same sequence
of \( d \)'s and \( N \)'s.

\begin{prop}
\label{porp:semicont}Let \( u_{1}(t),\dots ,u_{q}(t) \) be a family of unitaries
in \( M \), normalizing \( B\cong L^{\infty }[0,1] \), and for which \( \lim _{t\to 0}u_{j}(t)=u_{j}\in M \)
in the sense of \( * \)-strong topology. Then
\[
\chi (u_{1},\dots ,u_{n}:u_{n+1},\dots ,u_{q}\btw B)\geq \limsup _{t\to 0}\chi (u_{1}(t),\dots ,u_{n}(t):u_{n+1}(t),\dots ,u_{q}(t)\btw B).\]
The same conclusion holds for \( \chi ^{\omega } \).
\end{prop}
\begin{proof}
Let \( d_{1},\dots ,d_{m}\in B \) be fixed. Then because of Lemma \ref{lemma:sametop},
we have that, having fixed \( \epsilon ,\delta  \) and \( t_{0}>0 \), there
is a \( t<t_{0} \) and \( 0<\epsilon '<\epsilon  \), \( 0<\delta '<\delta  \)
for which
\[
\Gamma (u_{1}(t),\dots ,u_{q}(t):d_{1},\dots ,d_{n},\sigma ,\epsilon ',\delta ',l,d,N)\subset \Gamma (u_{1},\dots ,u_{q}:d_{1},\dots ,d_{n},\sigma ,\epsilon ,\delta ,l,d,N)\]
for all \( \sigma \in S_{N}^{n} \) and all \( d,N>0 \). The claimed inequality
now follows from the definition of \( \chi  \).
\end{proof}
\begin{prop}
\label{prop:finiteuwaut}Let \( \alpha  \) be an automorphism of \( B=L^{\infty }[0,1] \),
preserving Lebesgue measure. Let \( u \) be the unitary in \( B\rtimes _{\alpha }\mathbb {Z} \),
which implements \( \alpha  \). Let \( w \) independent of \( M \) and free
from \( u \) with amalgamation over \( M \). Then \( \chi (uw:w\btw M)\geq \chi (w) \)
and \( \chi (uw\btw M)\geq \chi (w) \).
\end{prop}
\begin{proof}
By Corollary \ref{corr:existanceofpermutation}, given \( d_{1},\dots ,d_{n}\in L^{\infty }[0,1] \),
\( \epsilon ,\delta ,l \), for \( N \) sufficiently large, there exists a
permutation \( \sigma \in S_{N} \), so that
\[
|d_{0}\sigma ^{g(1)}(d_{1}\sigma ^{g(1)}(\dots \sigma ^{g(k)}(d_{k})\dots )))-d_{0}\alpha ^{g(1)}(d_{1}\alpha ^{g(1)}(\dots \alpha ^{g(k)}(d_{k})\dots )))|_{\epsilon }<\delta ,\]
where \( \alpha =\Ad _{u} \). It follows that \( \sigma \cdot 1\in \Gamma ^{\sigma }(u:d_{1},\dots ,d_{m},\epsilon ,\delta ,l,d,N) \),
for all \( d \). Given \( \theta >0 \), for \( d \) large enough, there exists
a subset \( X\subset \Gamma (w;l,d,\delta )^{\oplus N} \), so that \( \mu (X)/\mu (\Gamma (w;l,d,\delta ))^{N}\geq 1-\theta  \),
and so that elements of \( X^{\oplus N} \) are free from \( \sigma  \) in
moments up to length \( l \) and degree \( \delta  \). It follows that the
set \( \{(\sigma \cdot x,x):x\in X\}\subset \Gamma ^{\sigma }(uw,w:d_{1},\dots ,d_{m},\epsilon ,\delta ,l,d,N) \).
The claimed inequality now follows from the definition of \( \chi  \).
\end{proof}
\begin{prop}
\label{prop:approxexists}Assume that \( u_{1},\dots ,u_{n}\in M \) normalize
\( B\cong L^{\infty }[0,1] \). Assume that \( u_{1},\dots ,u_{n},B \) have
f.d.a. Let \( w_{1},\dots ,w_{n} \) commute with \( B \), be independent from
\( B \), free with amalgamation over \( B \) from each other and free with
amalgamation over \( B \) from \( u_{1},\dots ,u_{n} \). Then
\[
\chi (w_{1}u_{1},\dots ,w_{n}u_{n}:w_{1},\dots ,w_{n}\btw B)\geq \sum _{j=1}^{n}\chi (w_{j}),\]

\[
\chi (w_{1}u_{1},\dots ,w_{n}u_{n}\btw B)\geq \sum _{j=1}^{n}\chi (w_{j}).\]

\end{prop}
\begin{proof}
Since \( u_{1},\dots ,u_{n} \) have f.d.a, given \( \epsilon ,\delta ,N_{0},d_{0} \),
there are \( d>d_{0} \), \( N>n_{0} \) so that for some \( \sigma  \) there
exists an element \( (U_{1},\dots ,U_{n})\in \Gamma ^{\sigma }(u_{1},\dots ,u_{n}:d_{1},\dots ,d_{n},\epsilon ,\delta ,l,d,N) \).
By the assumed freeness between \( w_{1},\dots ,w_{n} \) and \( u_{1},\dots ,u_{n} \),
we find that given \( \theta >0 \), for all \( N \) and \( d \) sufficiently
large, there is a subset \( \Gamma \subset \Gamma ^{\id }(w_{1},\dots ,w_{n}:d_{1},\dots ,d_{n},\epsilon ,\delta ,l,d,N) \),
so that \( \mu (\Gamma )/\mu (\Gamma ^{\id }(w_{1},\dots ,w_{n}:d_{1},\dots ,d_{n},\epsilon ,\delta ,l,d,N))\geq 1-\theta  \),
and so that
\[
(U_{1},\dots ,U_{n})\times \Gamma \subset \Gamma ^{\sigma \times \id }(u_{1},\dots ,u_{n},w_{1},\dots ,w_{n}:d_{1},\dots ,d_{n},\epsilon ,\delta ,l,d,N).\]
It follows that given \( \epsilon ',\delta ',l' \) there exist \( 0<\epsilon <\epsilon ',0<\delta <\delta ',l>l' \)
for which the image of the map
\[
\Gamma \ni (W_{1},\dots ,W_{n})\mapsto (W_{1}U_{1},\dots ,W_{n}U_{n})\]
lies in \( \Gamma ^{\sigma }(u_{1}w_{1},\dots ,u_{n}w_{n}:w_{1},\dots ,w_{n}:d_{1},\dots ,d_{m},\epsilon ,\delta ,l,d,N) \).
It follows after taking limits that
\begin{eqnarray*}
\chi (w_{1}u_{1},\dots ,w_{n}u_{n}\btw B) & \geq  & \chi (w_{1}u_{1},\dots ,w_{n}u_{n}:w_{1},\dots ,w_{n}\btw B)\\
 & \geq  & \chi (w_{1},\dots ,w_{n}\btw B).
\end{eqnarray*}
By the independence and freeness assumptions on \( w_{1},\dots ,w_{n} \) we
finally get
\[
\chi (w_{1},\dots ,w_{n}\btw B)=\sum \chi (w_{j}\btw B)=\sum \chi (w_{j}),\]
which is the desired estimate.
\end{proof}
\begin{prop}
\label{prop:changevar}Let \( u_{1},\dots ,u_{n},v_{1},\dots ,v_{m},w\in M \)
be in the normalizer of \( B \), and assume that \( y\in W^{*}(u_{1},\dots ,u_{n},B) \)
is a unitary, so that \( y \) normalizes \( B \). Then
\[
\chi (u_{1},\dots ,u_{n},w:v_{1},\dots ,v_{m}\btw B)=\chi (u_{1},\dots ,u_{n},yw:v_{1},\dots ,v_{m}\btw B).\]
The same statement holds for \( \chi  \) replaced by \( \chi ^{\omega } \).
The same conclusion holds even if \( m=\infty  \).
\end{prop}
The proof is only sketched, being for the most part exactly the same as the
proof of the change of variables formula (see \cite{dvv:entropy2}). Note that
in view of the assumption that \( u_{1},\dots ,u_{n} \) normalize \( B \),
one can approximate \( y \) by \( p(u_{1},\dots ,u_{n}) \), where \( p \)
is a polynomial with coefficients from \( B \) of the form \( p(t_{1},\dots ,t_{n})=\sum _{m}\sum f_{i_{1},\dots ,i_{m}}t_{i_{1}}\cdots t_{i_{m}} \),
with \( f_{\cdots }\in B \). It can be shown exactly as in \cite{dvv:entropy2}
that \( \chi (u_{1},\dots ,u_{n},w:v_{1},\dots ,v_{m})=\chi (u_{1},\dots ,u_{n},p(u_{1},\dots ,u_{n})w,v_{1},\dots ,v_{m}\btw B) \).
Taking limits gives by Proposition \ref{porp:semicont} the inequality \( \chi (u_{1},\dots ,u_{n},w:v_{1},\dots ,v_{m}\btw B)\leq \chi (u_{1},\dots ,u_{n},yw:v_{1},\dots ,v_{m}\btw B) \).
Replacing now \( y \) by \( y^{-1} \) gives the opposite inequality. The proof
in the case that \( m=\infty  \) is exactly the same.

\section{Free Dimension \protect\( \delta (\cdots :\cdots \btw B)\protect \).}

\begin{defn}
Given \( u_{1},\dots ,u_{n},v_{1},v_{2},\dots \in M \) normalizing \( L^{\infty }[0,1]\cong B\subset M \),
define
\[
\delta _{0}(u_{1},\dots ,u_{n}:v_{1},v_{2},\dots \btw B)=n-\liminf _{t\to 0}\frac{\chi (w_{1}(t)u_{1},\dots ,w_{n}(t)u_{n}:v_{1},v_{2},\dots ,w_{1}(t),\dots ,w_{n}(t)\btw B)}{\log t^{1/2}},\]
where \( w_{1}(t),\dots ,w_{n}(t) \) commute with \( B \), are independent
from \( B \), are free from each other with amalgamation over \( B \), and
are free from \( u_{1},\dots ,u_{n},v_{1},v_{2},\dots  \) with amalgamation
over \( B \), and are such that \( w_{j}(t) \) is \( * \)-distributed as
the multiplicative free Brownian motion started at identity and evaluated at
time \( t \). Here we \emph{allow} there to be an infinite set of \( v_{1},v_{2},\dots  \).

Define similarly \( \delta _{0}^{\omega } \) by replacing \( \chi  \) with
\( \chi ^{\omega } \). Finally, for an element \( \kappa \in \beta ((0,1])\setminus (0,1] \),
define \( \delta ^{\omega }_{0,\kappa }(u_{1},\dots ,u_{n}\btw B) \) by replacing
\( \liminf  \) in the definition of \( \delta  \) with \( \lim _{t\to \kappa } \).

Define also
\[
\delta (u_{1},\dots ,u_{n}:v_{1},v_{2}\dots \btw B)=n-\liminf _{t\to 0}\frac{\chi (w_{1}(t)u_{1},\dots ,w_{n}(t)u_{n}:v_{1},v_{2},\dots \btw B)}{\log t^{1/2}},\]
and \( \delta ^{\omega } \), \( \delta ^{\omega }_{k} \) in the obvious way.
\end{defn}
\begin{prop}
\label{prop:fdiminpresence}If \( w\in W^{*}(u_{1},\dots ,u_{n}) \), then \( \delta _{0}(u_{1},\dots ,u_{n}\btw B)=\delta _{0}(u_{1},\dots ,u_{n}:w\btw B) \).
\end{prop}
\begin{proof}
It is sufficient to prove that, with the same notation as in the definition
of \( \delta _{0} \), 
\begin{eqnarray}
\chi (u_{1}w_{1}(t),\dots ,u_{n}w_{n}(t):w_{1}(t),\dots ,w_{n}(t),y_{1},y_{2},\dots \btw B) &  & \label{eqn:changeofpresence} \\
=\chi (u_{1}w_{1}(t),\dots ,u_{n}w_{n}(t):w_{1}(t),\dots ,w_{n}(t),v,y_{1},y_{2},\dots \btw B) &  & 
\end{eqnarray}
(we caution the reader that the quantity on the left involves entropy in the
presence of an infinite number of variables). The inequality \( \leq  \) is
clear. To prove the opposite inequality, fix \( \delta >0 \), and choose \( r>0 \)
so that \( |E_{B}(|u-p(y_{1},\dots ,y_{r})|^{2}|_{\epsilon }<\delta  \) for
some non-commutative polynomial \( p \) with coefficients from \( B \). Then
one has the inclusion 
\begin{eqnarray*}
\Gamma (u_{1}w_{1}(t),\dots ,u_{n}w_{n}(t):w_{1}(t),\dots ,w_{n}(t),y_{1},\dots ,y_{q}:d_{1},\dots ,d_{m},\sigma ,\epsilon ,\delta ,l,d,N)\subset  &  & \\
\Gamma (u_{1}w_{1}(t),\dots ,u_{n}w_{n}(t):w_{1}(t),\dots ,w_{n}(t),y_{1},\dots ,y_{q},w:d_{1},\dots ,d_{m},\sigma ,l\epsilon ,2\delta ,l,d,N) &  & 
\end{eqnarray*}
for all \( q\geq r \). Taking limits gives the opposite inequality, and hence
implies (\ref{eqn:changeofpresence}). 
\end{proof}
\begin{prop}
\label{prop:stupidfdimbounds}\( \delta (u_{1},\dots ,u_{n}\btw B)\leq \sum \delta (u_{j}\btw B)\leq n \).
Moreover, if \( (u_{1},\dots ,u_{n},B) \) has f.d.a., then \( \delta (u_{1},\dots ,u_{n}\btw B)\geq 0 \).
In particular, for a single unitary \( u \) normalizing \( B \) we always
have \( 0\leq \delta (u\btw B)\leq 1 \). The same statements hold true for
\( \delta _{0} \), \( \delta _{0}^{\omega } \), \( \delta ^{\omega }_{0,\kappa } \),
\( \delta ^{\omega } \) and \( \delta ^{\omega }_{\kappa } \).
\end{prop}
\begin{proof}
The first inequality follows from \( \chi (v_{1},\dots ,v_{n}:w_{1},\dots ,w_{n}\btw B)\leq \sum \chi (v_{j}:w_{j}\btw B)\leq 0 \)
(note that \( \log t<0 \) for \( t \) close to zero). The second inequality
follows (under the assumptions of the hypothesis) from
\[
\chi (w_{1}(t)u_{1},\dots ,w_{n}(t)u_{n}:w_{1}(t),\dots ,w_{n}(t))\geq \sum _{j=1}^{n}\chi (w_{j}(t))=n\chi (w_{1}(t))\]
and from
\[
\lim _{t\to 0}\frac{\chi (w_{1}(t))}{\log t^{1/2}}=1.\]
(see \cite{shlyakht:cost}). 

The statement for one unitary follows from Corollary \ref{corr:existanceofpermutation}.
\end{proof}
\begin{rem}
It is easily seen that the condition \( \delta (u_{1},\dots ,u_{n}\btw B)\geq 0 \)
is equivalent to the assumption that \( (u_{1},\dots ,u_{n},B) \) has f.d.a.
(see Definition \ref{def:FDA}). Here \( \delta  \) can be replaced with \( \delta _{0} \),
\( \delta _{0}^{\omega } \), \( \delta ^{\omega }_{0,\kappa } \), \( \delta ^{\omega } \)
and \( \delta ^{\omega }_{\kappa } \).
\end{rem}
\begin{prop}
If the families \( (u_{1},\dots ,u_{n}),(v_{1},\dots ,v_{m}) \) are free with
amalgamation over \( B \), then
\[
\delta ^{\omega }_{\kappa }(u_{1},\dots ,u_{n},v_{1},\dots ,v_{n}\btw B)=\delta ^{\omega }_{k}(u_{1},\dots ,u_{n}\btw B)+\delta _{\kappa }^{\omega }(v_{1},\dots ,v_{m}\btw B).\]
 The same statement holds true for \( \delta _{0,\kappa }^{\omega } \).
\end{prop}
\begin{proof}
This follows from Proposition \ref{prop:freeadditivity}. 
\end{proof}
Note that the use of \( \lim _{t\to \kappa } \) in the definition of \( \delta ^{\omega }_{\kappa } \)
and \( \delta ^{\omega }_{0,\kappa } \) is crucial: otherwise, there is no
reason that additivity of free entropy \( \chi ^{\omega } \) translates into
additivity of free dimension, since we do not know if \( \liminf  \) in the
definition of free entropy is in general a limit.

\begin{prop}
\label{prop:indepofgenoneineq}Assume that \( u_{1},\dots ,u_{n}\in M \), \( v_{n+1},\dots ,v_{d}\in M \)
are unitaries normalizing \( D \). Let \( w_{1}(t),\dots ,w_{n}(t) \) be unitaries,
independent from \( B \), \( * \)-free with amalgamation over \( B \) from
each other and from \( u_{1},\dots ,u_{n},v_{n+1},\dots ,v_{d} \), and such
that each \( w_{j}(t) \) is \( * \)-distributed as multiplicative free Brownian
motion started at identity and evaluated at time \( t \). Assume that for a
fixed family of projections \( p_{n+1},\dots ,p_{d}\in A \) so that \( \tau (p_{j})=1-\rho _{j} \),
\( n<j\leq d \), and for each \( t>0 \) there exist unitaries \( P_{n+1}(t),\dots ,P_{d}(t)\in W^{*}(B,u_{1}w_{1}(t),\dots ,u_{n}w_{n}(t)) \),
so that:
\begin{enumerate}
\item \( P_{j}(t) \) normalizes \( B \);
\item \( P_{j}(t) \) commutes with \( p_{j}B \);
\item for all \( 0<s<1 \), \( \Vert E_{B}(|p_{j}P_{j}(t)v^{*}-p_{j}|^{2})\Vert ^{1/2}=O(t^{s/2}) \).
\end{enumerate}
Then
\[
\delta _{\kappa }^{\omega }(u_{1},\dots ,u_{n},v_{n+1},\dots ,v_{d}:y_{1},y_{2},\dots \btw B)\leq \delta _{k}^{\omega }(u_{1},\dots ,u_{n}:y_{1},y_{2},\dots \btw B)-\sum _{j=n+1}^{d}\rho _{j}.\]
The same statement holds for \( \delta _{0} \), \( \delta ^{\omega } \), \( \delta ^{\omega }_{0} \)
and \( \delta ^{\omega }_{0,\kappa } \). 
\end{prop}
\begin{proof}
It is sufficient to prove the statement for \( d=n+1 \). Write \( v \) for
\( v_{n+1} \), \( p \) for \( p_{n+1} \). Denote by \( P_{t} \) the unitary
\( P_{d}(t) \). The We have, using the definition of \( \delta  \), Proposition
\ref{prop:changevar} and subadditivity of entropy that
\begin{eqnarray*}
\delta (u_{1},\dots ,u_{n},v:y_{1},\dots \btw B) & = & n+1-\liminf _{t\to 0}\frac{\chi (u_{1}w_{1}(t),\dots ,u_{n}w_{n}(t),vw_{n+1}(t):y_{1},\dots \btw B)}{\log t^{1/2}}\\
 & = & n+1-\liminf _{t\to 0}\frac{\chi (u_{1}w_{1}(t),\dots ,u_{n}w_{n}(t),P_{t}^{*}vw_{n+1}(t):y_{1},\dots \btw B)}{\log t^{1/2}}\\
 & \leq  & n-\liminf _{t\to 0}\frac{\chi (u_{1}w_{1}(t),\dots ,u_{n}w_{n}(t):y_{1},\dots \btw B)}{\log t^{1/2}}\\
 &  & +1-\liminf _{t\to 0}\frac{\chi (P_{t}^{*}vw_{n+1}(t)\btw B)}{\log t^{1/2}}\\
 & = & \delta (u_{1},\dots ,u_{n}:y_{1},\dots \btw B)+1-\liminf _{t\to 0}\frac{\chi (P_{t}^{*}vw_{n+1}(t)\btw B)}{\log t^{1/2}}.
\end{eqnarray*}
Since \( P_{t}^{*}vw_{n+1}(t) \) commutes \( p \), we have by Proposition
\ref{prop:additivityundercompression}, Proposition \ref{Prop:independence}
and Proposition \ref{prop:estimateonchi} that for some constant \( D \) independent
of \( t \),
\begin{eqnarray*}
\chi (Q_{t}^{*}vw_{n+1}(t)\btw B) & = & \tau (p)\chi (pP_{t}^{*}vw_{n+1}(t)\btw pB)+(1-\tau (p))\chi ((1-p)Q_{t}^{*}vw_{n+1}(t)\btw (1-p)B)\\
 & \leq  & \tau (p)\chi (pP_{t}^{*}vw_{n+1}(t)\btw pB)\\
 & \leq  & \tau (p)\log t^{s/2}+\tau (p)D=\rho s\log t^{1/2}+O(\log t^{1/2}),
\end{eqnarray*}
since \( \Vert E_{B}(|pP_{t}^{*}vw_{n+1}(t)-p|^{2})\Vert ^{1/2}\leq \Vert E_{B}(|pP_{t}^{*}vw_{n+1}(t)-pw_{n+1}(t)|^{2})\Vert ^{1/2}+\Vert |p(w_{n+1}(t)-1)|\Vert =O(t^{s/2})+O(t^{1/2}) \).
It now follows that
\begin{eqnarray*}
\delta (u_{1},\dots ,u_{n},v:y_{1},\dots \btw B) & \leq  & \delta (u_{1},\dots ,u_{n}:y_{1},\dots \btw B)-\liminf _{t\to 0}\frac{\rho s\log t^{1/2}+O(\log t^{1/2})}{\log t^{1/2}}\\
 & = & \delta (u_{1},\dots ,u_{n}:y_{1},\dots \btw B)-\rho s.
\end{eqnarray*}
Since \( 0<s<1 \) was arbitrary, this implies the desired inequality. The proof
for \( \delta ^{\omega } \), \( \delta _{0} \), etc. is the same.
\end{proof}
\begin{prop}
\label{prop:indepofgenineq2}Assume that \( v_{1},\dots ,v_{m}\in W^{*}(u_{1},\dots ,u_{n},y_{1},y_{2},\dots ,B)\cap \mathcal{N}(B) \).
Then
\[
\delta _{0}(u_{1},\dots ,u_{n}:y_{1},y_{2},\dots \btw B)\leq \delta _{0}(u_{1},\dots ,u_{n},v_{1},\dots ,v_{m}:\btw B).\]
The same inequality is true for \( \delta ^{\omega }_{0} \) and \( \delta ^{\omega }_{0,\kappa } \).
In particular,
\[
\delta _{0}(u_{1},\dots ,u_{n}\btw B)\leq \delta _{0}(u_{1},\dots ,u_{n},v_{1},\dots ,v_{m}\btw B)\]
for all \( v_{1},\dots ,v_{m}\in W^{*}(B,u_{1},\dots ,u_{n})\cap \mathcal{N}(B) \).
\end{prop}
The proof is essentially identical to that of \cite[Theorem 4.3]{dvv:improvedrandom},
using Proposition \ref{prop:changeofvars} and Corollary \ref{corr:gammasalmostfree},
but we will provide it for completeness.

\begin{proof}
It is sufficient to prove the statement for \( m=1 \). Henceforth denote \( v_{1} \)
by \( v \). By Proposition \ref{prop:fdiminpresence}, we have that \( \delta _{0}(u_{1},\dots ,u_{n}:y_{1},y_{2}\btw B)=\delta _{0}(u_{1},\dots ,u_{n}:v,y_{1},y_{2},\dots \btw B) \).
Therefore, under the hypothesis of the Proposition, we have the inequality
\begin{eqnarray*}
\delta _{0}(u_{1},\dots ,u_{n}:y_{1},y_{2},\dots \btw B) & = & \delta _{0}(u_{1},\dots ,u_{n}:v,y_{1},y_{2},\dots \btw B)\\
 & \leq  & \delta _{0}(u_{1},\dots ,u_{n}:v\btw B).
\end{eqnarray*}
Thus, to conclude the proof, it is therefore sufficient to prove that \( \delta _{0}(u_{1},\dots ,u_{n}:v\btw B)\leq \delta _{0}(u_{1},\dots ,u_{n},v\btw B) \). 

Since \( \delta _{0}(u_{1},\dots ,u_{n},1:v\btw B)=\delta _{0}(u_{1},\dots ,u_{n}:v\btw B) \)
because \( 1 \) is free from \( u_{1},\dots ,u_{n},v \) with amalgamation
over \( B \), and \( \delta (1\btw B)=0 \), it follows that we must prove

\[
\delta _{0}(u_{1},\dots ,u_{n},1:v\btw B)\leq \delta _{0}(u_{1},\dots ,u_{n},v\btw B).\]
Thus it would be sufficient to prove the inequality
\begin{eqnarray*}
\chi (u_{1}w_{1}(t),\dots ,u_{n}w_{n}(t),w_{n+1}(t):w_{1}(t),\dots ,w_{n+1}(t),v\btw B)\leq  &  & \\
\chi (u_{1}w_{1}(t),\dots ,u_{n}(t)w_{n}(t),vw_{n+1}(t):w_{1}(t),\dots ,w_{n+1}(t)\btw B). &  & 
\end{eqnarray*}
Given \( \rho >0 \), there exists a Borel map \( G_{d,N} \), assuming a finite
number of values, from the set
\[
\Gamma (u_{1}w_{1}(t),\dots ,u_{n}w_{n}(t):w_{1}(t),\dots ,w_{n}(t),v:d_{1},\dots ,d_{m},\sigma ,\epsilon ,\delta ,l,d,N)\]
to the set 
\[
\Gamma (u_{1}w_{1}(t),\dots ,u_{n}w_{n}(t),w_{1}(t),\dots ,w_{n}(t),v:d_{1},\dots ,d_{m},\sigma ,\epsilon ,\delta ,l,d,N)\]
having the form
\begin{eqnarray*}
G_{d,N}(U_{1},\dots ,U_{n})=\big (f_{1}^{d,N}(U_{1},\dots ,U_{n}),\dots ,f_{n}^{d,N}(U_{1},\dots ,U_{n}), &  & \\
g_{1}^{d,N}(U_{1},\dots ,U_{n}),g_{n}^{d,N}(U_{1},\dots ,U_{n}), &  & \\
h^{d,N}(U_{1},\dots ,U_{n})\big ), &  & 
\end{eqnarray*}
so that \( |E_{\Delta _{N}}(|f_{k}^{d,N}(U_{1},\dots ,U_{n})-U_{k}|^{2})^{1/2}|_{\epsilon }\leq \rho  \)
for all \( 1\leq k\leq n \). 

Moreover, since \( w_{n+1}(t) \) is free with amalgamation over \( B \) from
\( u_{1},\dots ,u_{n},w_{1}(t),\dots ,w_{n}(t),v \), there exists a subset
\( \Omega (d,N) \) in \( \Gamma (u_{1}w_{1}(t),\dots ,u_{n}w_{n}(t),w_{n+1}(t):w_{1}(t),\dots ,w_{n}(t):d_{1},\dots ,d_{m},\sigma ,\epsilon ,\delta ,l,d,N)\times \Gamma (w_{n+1}(t):w_{1}(t),\dots ,w_{n}(t),d_{1},\dots ,d_{m},\sigma ,\epsilon ,\delta ,l,d,N) \),
so that
\[
\lim _{d}\frac{\mu (\Omega _{d,N})}{\begin{array}{c}
\mu (\Gamma (u_{1}w_{1}(t),\dots ,u_{n}w_{n}(t),w_{n+1}(t):w_{1}(t),\dots ,w_{n}(t),v:d_{1},\dots ,d_{m},\sigma ,\epsilon ,\delta ,l,d,N))\\
\Gamma (w_{n+1}(t):w_{1}(t),\dots ,w_{n}(t),d_{1},\dots ,d_{m},\sigma ,\epsilon ,\delta ,l,d,N))
\end{array}\times }=1\]
and so that for all \( U_{1},\dots ,U_{n},W\in \Omega  \), we have that
\[
(G(U_{1},\dots ,U_{n}),W)\in \Gamma (u_{1}w_{1}(t),\dots ,u_{n}w_{n}(t),w_{1}(t),\dots ,w_{n}(t),v,w_{n+1}(t):d_{1},\dots ,d_{m},\sigma ,\epsilon ,\delta ,l,d,N)\]

In particular, for \( U_{1},\dots ,U_{n},W\in \Omega  \), the values of the
map
\[
H(U_{1},\dots ,U_{n},W)=(U_{1},\dots ,U_{n},h^{d,N}(U_{1},\dots ,U_{n})W)\]
lie in the set
\[
\Gamma (u_{1}w_{1}(t),\dots ,u_{n}w_{n}(t),vw_{n+1}(t):w_{1}(t),\dots ,w_{n}(t),w_{n+1}(t):d_{1},\dots ,d_{m},\sigma ,l\epsilon ,\delta +\rho ,l,d,N).\]
Since this map preserves Haar measure on the unitary group, we conclude, after
passing to limits, that
\begin{eqnarray*}
\chi (u_{1}w_{1}(t),\dots ,u_{n}w_{n}(t),w_{n+1}(t):w_{1}(t),\dots ,w_{n+1}(t),v\btw B)\leq  &  & \\
\chi (u_{1}w_{1}(t),\dots ,u_{n}(t)w_{n}(t),vw_{n+1}(t):w_{1}(t),\dots ,w_{n+1}(t)\btw B), &  & 
\end{eqnarray*}
thus finishing the proof.
\end{proof}
\begin{defn}
Let \( u_{1},\dots ,u_{n},v_{1},\dots ,v_{m}\in M \) be in the normalizer of
\( B\cong L^{\infty }[0,1] \). We say that \( u_{1},\dots ,u_{n} \) and \( v_{1},\dots ,v_{m} \)
are orbit-equivalent, if there are projections \( p^{(j)}_{k} \), \( q^{(l)}_{s} \),
\( 1\leq j\leq n \), \( 1\leq s\leq m \), \( 1\leq k\leq N(j) \), \( 1\leq s\leq M(l) \)
with possibly \( N(j) \) or \( M(l)=\infty  \), words \( g^{(j)}_{k} \) consisting
of letters from \( u_{1},\dots ,u_{n},u_{1}^{*},\dots ,u_{n}^{*} \), and words
\( h_{k}^{(j)} \) consisting of letters from \( v_{1},\dots ,v_{m},v_{1}^{*},\dots ,v_{m}^{*} \),
so that
\[
v_{j}=\sum _{k=1}^{N(j)}p^{(j)}_{k}g_{k}^{(j)},\qquad u_{l}=\sum _{s=1}^{M(l)}q_{s}^{(l)}h_{s}^{(l)}.\]
(In particular, one must have \( \sum _{k}p_{k}^{(j)}=\sum _{l}q_{l}^{(j)}=1 \)).
\end{defn}
The results of Feldman and Moore \cite{feldman-moore} imply that if \( B\subset M \)
is a Cartan subalgebra, then \( u_{1},\dots ,u_{n} \) in the normalizer of
\( B \) are orbit-equivalent to \( v_{1},\dots ,v_{m} \) in the normalizer
of \( B \) iff \( W^{*}(u_{1},\dots ,u_{n},B)=W^{*}(v_{1},\dots ,v_{m},B) \).
This is the case, for example, if \( M=W^{*}(X,R) \) is the von Neumann algebra
of a measurable equivalence relation \( R \) on a measure space \( X \), and
\( B\subset M \) is the canonical copy of \( L^{\infty }(X) \) in \( W^{*}(X,R) \).

\begin{prop}
\label{prop:independenceofgenerators}Let \( u_{1},\dots ,u_{n}\in M \), \( v_{1},\dots ,v_{m}\in M \)
be unitaries normalizing \( B \). Assume that \( u_{1},\dots ,u_{n} \) and
\( v_{1},\dots ,v_{m} \) are orbit-equivalent over \( B \). Then \( \delta _{0}(u_{1},\dots ,u_{n}\btw B)=\delta _{0}(v_{1},\dots ,v_{m}\btw B) \).
The same conclusion holds for \( \delta _{0}^{\omega } \) and \( \delta ^{\omega }_{0,\kappa } \).
\end{prop}
\begin{proof}
Note that under the orbit-equivalence assumptions, for all \( 0<\rho _{j}<1 \),
\( j=1,\dots ,m \), there exist polynomials \( P_{j} \) with coefficients
from \( B \) having the form \( P_{j}(z_{1},\dots ,z_{n})=\sum _{k=1}^{N_{j}}q_{k}^{(j)}z_{i_{1}}^{\pm 1}\dots z_{i_{t(k)}}^{\pm 1} \),
where \( q_{k}(j) \) are orthogonal projections, so that \( p_{j}v_{j}=p_{j}P_{j}(u_{1},\dots ,u_{n}) \),
where \( p_{j}=\sum _{k=1}^{N_{j}}q_{k}^{(j)} \) and \( \tau (p_{j})=1-\rho _{j} \).
In particular, \( p_{j}P_{j}(u_{1}w_{1}(t),\dots ,w_{n}(t))v_{j}^{*} \) commutes
with \( q_{k}^{j}B \) whenever \( w_{1}(t),\dots ,w_{n}(t) \) are unitaries
and commute with \( B \). Take \( w_{1}(t),\dots ,w_{n}(t) \) to be free Brownian
motion, as in the definition of the free dimension \( \delta (\cdot \btw B) \).
Since \( u_{i}w_{i}(t) \), \( i=1,\dots ,n \) normalize \( B \) and define
the same automorphisms of \( B \) as \( u_{1},\dots ,u_{n} \), it follows
that there exists unitaries \( P_{j}(t)\in W^{*}(u_{1}w_{1}(t),\dots ,u_{n}w_{n}(t)) \),
\( j=1,\dots ,m \), normalizing \( B \), so that \( p_{j}P_{j}(t)=p_{j}P_{j}(u_{1}w_{1}(t),\dots ,u_{n}w_{n}(t)) \)
(one can simply choose any extension of the isometry \( p_{j}P_{j}(u_{1}w_{1}(t),\dots ,u_{n}w_{n}(t))\in W^{*}(u_{1}w_{1}(t),\dots ,u_{n}w_{n}(t)) \)
to a unitary normalizing \( B \)). Therefore, since \( \Vert w_{j}(t)-1\Vert =O(t^{1/2}) \)
(cf. \cite{biane:brownian}),
\[
\Vert p_{j}P_{j}(t)v_{j}^{*}-p_{j}\Vert =O(t^{1/2}),\]
hence \( {\Vert E_{B}(|p_{j}P_{j}(t)v_{j}^{*}-p_{j}|^{2})\Vert }^{1/2}=O(t^{1/2}) \).
It follows that the hypothesis of Proposition \ref{prop:indepofgenoneineq}
is satisfied, and hence \( \delta _{0}(u_{1},\dots ,u_{n},v_{1},\dots ,v_{m}\btw B)\leq \delta _{0}(u_{1},\dots ,u_{n}\btw B)-\sum \rho _{j} \).
By Proposition, \ref{prop:indepofgenineq2} we get also that \( \delta _{0}(u_{1},\dots ,u_{n})\leq \delta _{0}(u_{1},\dots ,u_{n},v_{1},\dots ,v_{m}\btw B) \).
Since \( \rho _{j} \) are arbitrary, we get that \( \delta _{0}(u_{1},\dots ,u_{n}\btw B)=\delta _{0}(u_{1},\dots ,u_{n},v_{1},\dots ,v_{m}\btw B) \).
Reversing the roles of \( u_{1},\dots ,u_{n} \) and \( v_{1},\dots ,v_{m} \)
gives finally that
\[
\delta _{0}(u_{1},\dots ,u_{n}\btw B)=\delta _{0}(u_{1},\dots ,u_{n},v_{1},\dots ,v_{m}\btw B)=\delta _{0}(v_{1},\dots ,v_{m}\btw B).\]

\end{proof}

\section{Computation of \protect\( \delta \protect \) for certain variables.}

\begin{lem}
\label{lemma:aboutcyclic}Let \( v(nt) \) be a unitary, classically independent
from an algebra \( A \) with a trace \( \tau  \). Let \( n \) be a positive
integer, and consider \( M=A\otimes M_{n\times n} \), with the trace \( \tau \otimes \frac{1}{n}\Tr  \).
Assume that \( v(nt) \) is \( * \)-distributed as a multiplicative free Brownian
motion started at identity and evaluated at time \( nt \). Let \( u_{1},\dots ,u_{n-1},u_{n} \)
be Haar unitaries, which are classically \( * \)-independent from \( A \)
and free from each other over \( A \). Let \( w_{1}=u_{1},\dots ,\, \, w_{n-1}=u_{n-1} \)
and \( w_{n}=(w_{1}\cdots w_{n-1})^{-1}v(nt) \), Consider the unitary
\[
Y(t)=\left( \begin{array}{ccccc}
0 & w_{1} & 0 & \cdots  & 0\\
0 & 0 & w_{2} & 0 & \vdots \\
\vdots  & \vdots  & \ddots  & \ddots  & \vdots \\
0 & 0 & \cdots  & 0 & w_{n-1}\\
w_{n} & 0 & \cdots  & 0 & 0
\end{array}\right) \in M.\]
Let \( B\cong A\otimes \mathbb {C}^{n} \) be the algebra of diagonal matrices
in \( M \) with entries from \( A \). Consider the automorphism of \( B \)
given by \( \id \otimes \sigma  \), where \( \sigma  \) is the cyclic permutation
on \( \mathbb {C}^{n} \). Let \( \sigma \in M=B\rtimes _{\id \otimes \sigma }\mathbb {Z}_{n} \)
be the canonical unitary implementing \( \id \otimes \sigma  \), and let \( w(t) \)
be a unitary, independent of \( B \), free from \( B\rtimes _{\id \otimes \sigma }\mathbb {Z}_{n} \)
with amalgamation over \( B \), and \( * \)-distributed as the free Brownian
motion started at identity and evaluated at time \( t \). 

Then the \( B \)-valued distribution of \( Y(t) \) is the same as the \( B \)-valued
distribution of \( u\sigma (t) \). Moreover, \( Y(t) \) is free from \( M \)
with amalgamation over \( B \).
\end{lem}
\begin{proof}
Note that the unitary
\[
U=\left( \begin{array}{ccc}
u_{1} & \cdots  & 0\\
\vdots  & \ddots  & \vdots \\
0 & \cdots  & u_{n}
\end{array}\right) \]
is free from \( B\rtimes _{\id \otimes \sigma }\mathbb {Z}_{n}\cong M \) with
amalgamation over \( B \), and is independent from \( B \). In our identification
of \( B\rtimes _{\id \otimes \sigma }\mathbb {Z}_{n} \) with \( A\otimes M_{n\times n} \)
the unitary \( u \) is identified with the matrix
\[
\Sigma =\left( \begin{array}{cccc}
0 & 1 & \cdots  & 0\\
0 & \ddots  & \ddots  & \vdots \\
\vdots  & \ddots  & \ddots  & 1\\
1 & \cdots  & 0 & 0
\end{array}\right) .\]
Lastly, if \( w_{1}(t),\dots ,w_{n}(t) \) are each \( * \)-distributed as
\( w(t) \), are independent from \( A \) and are \( * \)-free over \( A \),
then the matrix
\[
W(t)=\left( \begin{array}{ccc}
w_{1}(t) & \cdots  & 0\\
\vdots  & \ddots  & \vdots \\
0 & \cdots  & w_{n}(t)
\end{array}\right) \]
is independent from \( B \), is \( * \)-distributed in the same way as \( w(t) \)
and is \( * \)-free with amalgamation over \( B \) from the \( * \)-algebra
generated by \( U \) and \( \Sigma  \). It follows that the \( B \)-valued
distribution of \( \Sigma W(t) \) is the same as the \( B \)-valued distribution
of \( \sigma w(t) \). Since \( U \) is free from \( \Sigma W(t) \) over \( B \),
and because \( U \) is a Haar unitary, independent from \( B \), it follows
that \( U\Sigma W(t)U^{*} \) is free from \( M \) over \( B \), and has the
same \( B \)-valued distribution as \( \sigma w(t) \).

Write \( Z(t)=U\Sigma W(t)U^{*} \). It remains to show that \( Y(t) \) and
\( Z(t) \) have the same \( M \)-valued \( * \)-distributions. Indeed, that
would imply that the \( B \)-valued distribution of \( Y(t) \) is the same
as that of \( Z(t) \) (hence the same as \( \sigma w(t) \)), and also that
\( Y(t) \) is \( * \)-free from \( M \) over \( B \), since \( Z(t) \)
is \( * \)-free from \( M \) over \( B \). As a matrix,
\[
Z(t)=\left( \begin{array}{ccccc}
0 & u_{1}w_{1}(t)u_{2}^{*} & 0 & \cdots  & 0\\
0 & 0 & u_{2}w_{2}(t)u_{3}^{*} & 0 & \vdots \\
\vdots  & \vdots  & \ddots  & \ddots  & \vdots \\
0 & 0 & \cdots  & 0 & u_{n-1}w_{n-1}(t)u_{n}^{*}\\
u_{n}w_{n}(t)u_{1}^{*} & 0 & \cdots  & 0 & 0
\end{array}\right) .\]

To prove that the \( M \)-valued \( * \)-distributions of \( Y(t) \) and
\( Z(t) \) are the same, it is sufficient to prove that the families of their
entries have the same joint \( * \)-distributions; i.e., that the joint \( * \)-distribution
of family \( (w_{1},\dots ,w_{n}) \) is the same as that of \( (u_{1}w_{1}(t)u_{2}^{*},u_{2}w_{2}(t)u_{3}^{*},\dots ,u_{n}w_{n}(t)u_{1}^{*}) \).
Write \( z_{1}=u_{1}w_{1}(t)u_{2}^{*},\dots ,z_{n}=u_{n}w_{n}(t)u_{1}^{*} \).
Hence it is sufficient to prove that: (i) \( z_{1},\dots ,z_{n-1} \) are Haar
unitaries, independent from \( A \) and \( * \)-free with amalgamation over
\( A \); (ii) \( v=z_{1}\cdots z_{n} \) is \( * \)-free from \( z_{1},\dots ,z_{n-1} \)
over \( A \) and (iii) \( v \) has the same \( A \)-valued \( * \)-distribution
as \( v(nt) \).

To prove (i), notice that we can, by replacing each \( w_{j}(t) \) with \( r_{j}w_{j}(t)r_{j}^{-1} \)
where \( r_{1},\dots ,r_{n} \) are Haar unitaries, independent from \( A \)
and \( * \)-free from each other and from \( w_{1}(t),\dots ,w_{n}(t),u_{1},\dots ,u_{n} \)
with amalgamation over \( A \), without changing the joint \( A \)-valued
\( * \)-distribution of the family, replace \( (z_{1},\dots ,z_{n-1}) \) by
\( (u_{1}r_{1}w_{1}(t)(u_{2}r_{1})^{*},\dots ,u_{n-1}r_{n-1}w_{n-1}(t)(u_{n}r_{n-1})^{*} \).
Since the unitaries \( (u_{1}r_{1},u_{2}r_{1},u_{2}r_{2},u_{2}r_{3},\dots ,u_{n-1}r_{n-1},u_{n}r_{n-1},w_{1}(t),\dots ,w_{n-1}(t)) \)
are \( * \)-free over \( A \), it follows that \( z_{1},\dots ,z_{n-1} \)
are \( * \)-free over \( A \). Clearly, each \( z_{j} \) is independent from
\( A \); and each \( z_{j} \) is a Haar unitary (note that we can always replace,
say, \( u_{j} \) by \( \exp (it)u_{j} \) for arbitrary \( t \), without changing
the joint distribution of \( z_{j} \)).

For the second claim, we have \( v=u_{1}w_{1}(t)\cdots w_{n}(t)u_{1}^{*} \).
Notice that \( u_{1} \) is \( * \)-free over \( A \) from \( (w_{1}(t)u_{2},\dots ,u_{n-1}w_{n-1}(t),w_{1}(t),\dots ,w_{n-1}(t)) \)
(which are all \( * \)-free over \( A \) among each other) and hence from
\( (u_{1}w_{1}(t)u_{2},\dots ,u_{n-1}w_{n-1}(t)u_{1}^{*},w_{1}(t),\dots ,w_{n-1}(t)) \).
Hence \( v \) is \( * \)-free over \( A \) from \( z_{1},\dots ,z_{n-1} \).

Lastly, \( v \) is clearly independent from \( A \), and has the same \( * \)-distribution
as \( w_{1}(t)\cdots w_{n}(t) \). Since \( w_{j}(t) \) are \( * \)-free and
form a multiplicative free Brownian motion, the \( * \)-distribution of \( w_{1}(t)\cdots w_{n}(t) \)
is the same as that of \( v(nt) \).
\end{proof}
\begin{prop}
\label{prop:fdimcyclic}Let \( n\in \mathbb {N} \) be fixed, and let \( \alpha  \)
be a free action of \( \mathbb {Z}_{n} \) on \( [0,1] \), and denote by \( u\in M=L^{\infty }[0,1]\rtimes _{\alpha }\mathbb {Z}_{n} \)
the associated unitary, implementing this action. Denote the canonical copy
of \( L^{\infty }[0,1]\subset M \) by \( B \). Then
\[
\delta ^{\omega }_{\kappa }(u\btw B)=1-\frac{1}{n},\]
independent of the choice of \( \omega  \) and \( \kappa  \); the same conclusion
holds for \( \delta  \) and \( \delta ^{\omega } \).

The same conclusion holds for \( \delta _{0} \), \( \delta ^{\omega }_{0} \)
and \( \delta ^{\omega }_{0,\kappa } \).
\end{prop}
\begin{proof}
We first prove the statement for \( \delta  \). We must prove that
\[
\lim _{t\to 0}\frac{\chi (w(t)u\btw B)}{\frac{1}{2}\log t}=\frac{1}{n}.\]

We shall prove that \( \chi (w(t)u\btw B)=\frac{1}{n}\chi (w(nt)) \), which
is sufficient, since
\[
2\lim _{t\to 0}\frac{\chi (w(nt))}{\log t}=2\lim _{r\to 0}\frac{\chi (w(r))}{\log r-\log n}=2\lim _{r\to 0}\frac{\chi (w(r))}{\log r}=1.\]

Choose cross-sections for the action of \( \mathbb {Z}_{n} \) on \( B \),
so that \( B\cong A\otimes \mathbb {C}^{n} \) and the action \( \alpha  \)
has the form \( \id \otimes \sigma  \) for a cyclic permutation \( \sigma  \)
or order \( n \) acting on \( \mathbb {C}^{n} \). Note that \( M\cong A\otimes M_{n\times n} \)
in such a way that identifies \( B \) with diagonal matrices in \( M \) with
values from \( A \), and \( u \) with the permutation matrix \( \sigma \in M_{n\times n} \).
Let \( v(nt),u_{1},\dots ,u_{n-1} \) be unitaries, independent from \( A \),
and free from each other over \( A \), and so that each \( u_{j} \) is a Haar
unitary, and \( v(nt) \) is \( * \)-distributed as a free multiplicative Brownian
motion started at identity and evaluated at time \( nt \). Let \( d_{1},\dots ,d_{r}\in A \)
be fixed, and let
\[
(V,U_{1},\dots ,U_{n-1})\in \Gamma (v(nt),u_{1}\dots ,u_{n-1}:d_{1}',\dots ,d_{r'}',\id ,\epsilon ',\delta ',l',d,N').\]
Set \( W_{1}=U_{1},\dots ,W_{n-1}=U_{n-1} \) and \( W_{n}=W_{1}\cdots W_{n-1}V \).
Let \( N \), \( \epsilon  \), \( \delta  \), \( d \) be given. For \( N \)
sufficiently large, we can write \( N=nN'+k \), where \( k<n \) and \( \frac{k}{N}<\frac{\epsilon }{2} \).
Then there exist \( \delta ' \), \( l' \), \( \epsilon ' \), \( r' \) and
\( d_{1}',\dots ,d_{r'}' \) for which the map
\begin{eqnarray*}
\Psi _{u}:\Gamma (v(nt),u_{1},\dots ,u_{n-1}:d_{1}',\dots ,d_{r'}',\id ,\epsilon ',\delta ',l',d,N')\ni (V,U_{1},\dots ,U_{n-1}) &  & \\
\mapsto \left( \begin{array}{ccccc}
0 & W_{1} & 0 & \cdots  & 0\\
\vdots  & \ddots  & \ddots  & \ddots  & \vdots \\
0 & \cdots  & 0 & W_{n-1} & 0\\
W_{n} & 0 & \cdots  & 0 & 0\\
0 & \cdots  & 0 & 0 & u
\end{array}\right) \in M_{nN'd+kd\times nN'd+kd} &  & 
\end{eqnarray*}
for a chosen matrix \( u\in M_{k} \) has values in \( \Gamma ^{\sigma \oplus \id _{k}}(w(t)u:d_{1},\dots ,d_{n},\epsilon ,\delta ,l,d,N) \),
and is injective. The union of its images over possible different \( u \) has
the same volume as \( \Gamma (v(nt),u_{1},\dots ,u_{n-1}:d_{1}',\dots ,d_{r'}',\id ,\epsilon ',\delta ',l',d,N')\times U(kd) \),
and is a subset of \( \Gamma ^{\sigma \oplus \id _{k}}(w(t)u:d_{1},\dots ,d_{n},\epsilon ,\delta ,l,d,N) \).
It follows that 
\[
\chi (w(t)u:d_{1},\dots ,d_{m},\epsilon ,\delta ,l,d,N)\geq \frac{N'}{nN'+k}\chi (v(nt),u_{1},\dots ,u_{n-1}:d_{1}',\dots ,d_{r'}',\epsilon ',\delta ',l',d,N'),\]
from which, after taking limits we get
\begin{equation}
\label{eqn:cyclicentropyinequality}
\chi (w(t)u\btw B)\geq \frac{1}{n}\chi (v(nt),u_{1},\dots ,u_{n-1}\btw A).
\end{equation}

Consider now the set \( \Gamma (w(t)u:d_{1},\dots ,d_{m},\bar{\sigma },\epsilon ,\delta ,l,d,N) \).
We may assume, for \( N \) large enough, that \( N=nN'+k \), \( k<n \), and
\( \bar{\sigma } \) has the form \( \sigma \oplus \id _{k} \), where \( \sigma \in M_{nN'\times nN'}\cong M_{n\times n}\otimes M_{N'\times N'} \)
has the form \( \id \otimes \sigma _{n} \), with \( \sigma _{n} \) a cyclic
permutation of order \( n \). Let \( \rho >0 \) be given. We may furthermore
assume by Lemma \ref{lemma:replacesigma} that for this choice of \( \bar{\sigma } \),
there exist \( \epsilon ''<\epsilon ,l''>l,\delta ''<\delta  \) for which \( \frac{1}{dN^{2}}\log \mu \Gamma (w(t)u:d_{1},\dots ,d_{m},\sigma ,\epsilon '',\delta '',l'',d,N) \)
is within \( \rho  \) of \( \chi (w(t)u:d_{1},\dots ,d_{m},\epsilon ,\delta ,l,d,N) \). 

Note that each element of \( \Gamma ^{\bar{\sigma }}(w(t)u:d_{1},\dots ,d_{n},\epsilon '',\delta '',l'',d,N) \)
lies in \( M_{nN'+k\times nN'+k}\otimes M_{d\times d} \) and can be represented
as a matrix
\begin{equation}
\label{eqn:U}
U=\left( \begin{array}{ccccc}
0 & W_{1} & 0 & \cdots  & 0\\
\vdots  & \ddots  & \ddots  & \ddots  & \vdots \\
0 & \cdots  & 0 & W_{n-1} & 0\\
W_{n} & 0 & \cdots  & 0 & 0\\
0 & \cdots  & 0 & 0 & u
\end{array}\right) 
\end{equation}
in which \( u\in M_{kd\times kd} \) and each \( W_{j}\in M_{dN'} \). By Corollary
\ref{corr:gammasalmostfree}, for large enough \( d \), there is a subset \( \bar{\Gamma } \)
of \( \Gamma =\Gamma (w(t)u:d_{1},\dots ,d_{n},\bar{\sigma },\epsilon '',\delta '',l'',d,N) \),
with \( \mu (\bar{\Gamma })/\mu (\Gamma )>\exp (-\rho ) \), so that each \( U\in \bar{\Gamma } \)
is \( \delta '',l'' \) free from the algebra \( M_{nN'+k\times nN'+k}\otimes 1 \),
and in particular, from \( \bar{\sigma }\oplus \id _{k} \). It follows that
given \( d_{1}',\dots ,d_{r'}'\in A \), \( \epsilon  \), \( \delta  \), \( l \),
there exist \( d_{1},\dots ,d_{r} \), \( \epsilon '<\epsilon '' \), \( l'>l'' \)
and \( \delta '<\delta '' \), so that if the matrix \( U \) above lies in
\( \bar{\Gamma }\cdot (\bar{\sigma }\oplus k) \), then \( (W_{1},\dots ,W_{n})\in \Gamma (v(nt),u_{1},\dots ,u_{n-1}:d_{1},\dots ,d_{r},\epsilon ',\delta ',l',d,N') \).
It follows that
\begin{eqnarray*}
\chi (w(t)u:d_{1},\dots ,d_{m},\epsilon ,\delta ,l,d,N)-\log (1-2\rho )\leq  &  & \\
\frac{N'}{nN'+k}\chi (v(nt),u_{1},\dots ,u_{n-1}:d_{1}',\dots ,d_{r'}',\epsilon ',\delta ',l',d,N'). &  & 
\end{eqnarray*}
Hence
\[
\chi (w(t)u\btw B)\leq \frac{1}{n}\chi (v(nt),u_{1},\dots ,u_{n-1}\btw A)+\log (1-2\rho ),\]
for \( \rho >0 \) arbitrarily small. Combining this with (\ref{eqn:cyclicentropyinequality})
gives, in view of independence and freeness assumptions:
\begin{eqnarray*}
\chi (w(t)u\btw B) & = & \frac{1}{n}\chi (v(nt),u_{1},\dots ,u_{n-1}\btw A),\\
 & = & \frac{1}{n}\chi (v(nt),u_{1},\dots ,u_{n-1})\\
 & = & \frac{1}{n}\chi (v(nt))+\frac{1}{n}\sum _{j=1}^{n-1}\chi (u_{j})\\
 & = & \frac{1}{n}\chi (v(nt))=\frac{1}{n}\chi (w(nt)),
\end{eqnarray*}
as we claimed.

The same proof can be modified to work for \( \delta _{0} \) instead; we point
out the necessary changes. We claim that
\[
\chi (w(t)u:w(t)\btw B)=\frac{1}{n}\chi (w_{1}(t)u_{1},\dots ,w_{n-1}(t)u_{n-1},w_{n}(t)u_{n-1}^{*}\cdots u_{1}^{*}:w_{1}(t),\dots ,w_{n}(t)),\]
where \( u_{1},\dots ,u_{n-1} \) are \( * \)-free Haar unitaries and \( w_{1}(t),\dots ,w_{n}(t) \)
are \( * \)-free unitaries, \( * \)-free from \( u_{1},\dots ,u_{n-1} \),
and each \( w_{j}(t) \) has the same distribution as free multiplicative Brownian
motion started from \( 1 \) and evaluated at time \( t \). The map \( \rho _{u} \)
which sends
\begin{eqnarray*}
(V_{1},\dots ,V_{n},W_{1},\dots ,W_{n})\in  &  & \\
\Gamma ^{\id }(w_{1}(t)u_{1},\dots ,w_{n-1}(t)u_{n-1},w_{n}(t)u_{n-1}^{*}\cdots u_{1}^{*},w_{1}(t),\dots ,w_{n}(t):d_{1}',\dots ,d_{r'}',\epsilon ,\delta d,l,N') &  & 
\end{eqnarray*}
to the pair of matrices
\[
\left( \left( \begin{array}{ccccc}
0 & V_{1} & 0 & \cdots  & 0\\
\vdots  & \ddots  & \ddots  & \ddots  & \vdots \\
0 & \cdots  & 0 & V_{n-1} & 0\\
V_{n} & 0 & \cdots  & 0 & 0\\
0 & \cdots  & 0 & 0 & u
\end{array}\right) ,\left( \begin{array}{ccccc}
W_{1} & 0 & \cdots  & 0 & 0\\
0 & \ddots  & \ddots  & \vdots  & \vdots \\
\vdots  & \ddots  & \ddots  & 0 & 0\\
0 & \cdots  & 0 & W_{n} & 0\\
0 & \cdots  & 0 & 0 & 1_{kd}
\end{array}\right) \right) \]
has values in
\[
\Gamma ^{\sigma }(w(t)u,w(t):d_{1},\dots ,d_{r},\epsilon ,\delta ,l,nN'+k).\]
This gives, just like in the first part of the proof, the inequality 
\[
\chi (w(t)u:w(t)\btw B)\geq \frac{1}{n}\chi (w_{1}(t)u_{1},\dots ,w_{n-1}(t)u_{n-1},w_{n}(t)u_{n-1}^{*}\cdots u_{1}^{*}:w_{1}(t),\dots ,w_{n}(t)\btw A).\]
Conversely, we can assume that there is a subset \( \bar{\Gamma } \) of \( \Gamma ^{\sigma }(w(t)u:w(t):d_{1},\dots ,d_{r},\epsilon ,\delta ,l,nN'+k) \),
so that
\[
\mu (\bar{\Gamma })/\mu \Gamma ^{\sigma }(w(t)u:w(t):d_{1},\dots ,d_{r},\epsilon ,\delta ,l,nN'+k)\geq \exp (-\rho ),\]
and so that for all \( U\in \bar{\Gamma } \) there exists a matrix \( V \),
commuting with \( \Delta _{N} \), for which \( (U,V)\in \Gamma ^{\sigma }(w(t)u,w(t):d_{1},\dots ,d_{r},\epsilon ,\delta ,l,nN'+k) \)
and \( (U,V) \) are \( l,\delta  \)-free from \( \sigma  \) with amalgamation
over \( \Delta _{N} \). Then the map sending such a pair \( (U,V) \),
\[
U=\left( \begin{array}{ccccc}
0 & V_{1} & 0 & \cdots  & 0\\
\vdots  & \ddots  & \ddots  & \ddots  & \vdots \\
0 & \cdots  & 0 & V_{n-1} & 0\\
V_{n} & 0 & \cdots  & 0 & 0\\
0 & \cdots  & 0 & 0 & u
\end{array}\right) ,\qquad V=\left( \begin{array}{cccc}
W_{1} &  &  & \\
 & \ddots  &  & \\
 &  & W_{n} & \\
 &  &  & w
\end{array}\right) \]
to \( (V_{1},\dots ,V_{n},W_{1},\dots ,W_{n}) \) is valued in
\[
\Gamma ^{\id }(w_{1}(t)u_{1},\dots ,w_{n-1}(t)u_{n-1},w_{n}(t)u_{n-1}^{*}\cdots u_{1}^{*},w_{1}(t),\dots ,w_{n}(t):d_{1}',\dots ,d_{r'}',\epsilon ,\delta d,l,N').\]
To see this, observe that the family \( (w_{1}(t)u_{1},\dots ,w_{n-1}(t)u_{n-1},w_{n}(t)u_{n-1}^{*}\cdots u_{1}^{*},w_{1}(t),\dots ,w_{n}(t)) \)
has the same joint \( A \)-valued \( * \)-distribution as
\[
(u_{1}w_{1}(t)u_{2}^{*},\dots ,u_{n-1}w_{n-1}(t)u_{n}^{*},u_{n}w_{n}(t)u_{1}^{*},u_{1}w_{1}(t)u_{1}^{*},\dots ,u_{n}w_{n}(t)u_{n}^{*}).\]
Next, observe that (as in the proof of Lemma \ref{lemma:aboutcyclic}) that
the family of matrices
\[
\left( \begin{array}{cccc}
0 & u_{1}w_{1}(t)u_{2}^{*} & \cdots  & 0\\
\vdots  & \ddots  & \ddots  & \vdots \\
0 & \cdots  & 0 & u_{n-1}w_{n-1}(t)u_{n}^{*}\\
u_{n}w_{n}(t)u_{1}^{*} & \cdots  & 0 & 0
\end{array}\right) ,\left( \begin{array}{cccc}
u_{1}w_{1}(t)u_{1}^{*} &  &  & \\
 & \ddots  &  & \\
 &  & \ddots  & \\
 &  &  & u_{n}w_{n}(t)u_{n}^{*}
\end{array}\right) \]
are \( * \)-free with amalgamation over \( B \) from the permutation matrix
\[
\sigma =\left( \begin{array}{cccc}
0 & 1 & \cdots  & 0\\
\vdots  & \ddots  & \ddots  & \vdots \\
0 & \cdots  & 0 & 1\\
1 & \cdots  & 0 & 0
\end{array}\right) ,\]
since they can be written as \( UW\sigma U^{*} \) and \( U\sigma U^{*} \),
where
\[
U=\left( \begin{array}{ccc}
u_{1} &  & \\
 & \ddots  & \\
 &  & u_{n}
\end{array}\right) ,\qquad W=\left( \begin{array}{ccc}
w_{1}(t) &  & \\
 & \ddots  & \\
 &  & w_{n}(t)
\end{array}\right) .\]

From this it follows that \( \bar{\Gamma } \) can be embedded into
\[
\Gamma ^{\id }(w_{1}(t)u_{1},\dots ,w_{n-1}(t)u_{n-1},w_{n}(t)u_{n-1}^{*}\cdots u_{1}^{*}:w_{1}(t),\dots ,w_{n}(t):d_{1}',\dots ,d_{r'}',\epsilon ,\delta d,l,N');\]
arguing as in the first part of the proof now gives
\[
\chi (w(t)u:w(t)\btw B)\geq \frac{1}{n}\chi (w_{1}(t)u_{1},\dots ,w_{n-1}(t)u_{n-1},w_{n}(t)u_{n-1}^{*}\cdots u_{1}^{*}:w_{1}(t),\dots ,w_{n}(t)\btw A)\]
and hence,
\[
\chi (w(t)u:w(t)\btw B)=\frac{1}{n}\chi (w_{1}(t)u_{1},\dots ,w_{n-1}(t)u_{n-1},w_{n}(t)u_{n-1}^{*}\cdots u_{1}^{*}:w_{1}(t),\dots ,w_{n}(t)\btw A).\]
To finish the proof, we must compute
\begin{eqnarray*}
1-\frac{1}{n}\liminf _{t\to 0}\frac{\chi (w_{1}(t)u_{1},\dots ,w_{n-1}(t)u_{n-1},w_{n}(t)u_{n-1}^{*}\cdots u_{1}^{*}:w_{1}(t),\dots ,w_{n}(t))\btw A}{\log t^{1/2}} & = & \\
1+\frac{1}{n}(\delta _{0}(u_{1},\dots ,u_{n-1},u_{n-1}^{*},\dots ,u_{1}^{*}\btw A)-n) & = & \\
1+\frac{1}{n}(\delta _{0}(u_{1},\dots ,u_{n-1}\btw A)-n) & = & \\
1+\frac{1}{n}(\sum _{j=1}^{n-1}\delta _{0}(u_{j}\btw A)-n) & = & \\
1+\frac{1}{n}((n-1)-n)=1-\frac{1}{n}, &  & 
\end{eqnarray*}
where we use freeness of \( u_{1},\dots ,u_{n-1} \) with amalgamation over
\( A \), and Proposition \ref{prop:independenceofgenerators}.
\end{proof}

\section{Free Dimension of an Equivalence Relation and Cost.}

\begin{prop}
Let \( R \) be a measurable equivalence relation on a finite measure space
\( X \). Assume that \( R \) has a finite graphing; i.e., there are automorphisms
\( \alpha _{1},\dots ,\alpha _{n},\dots  \) of \( X \), which generate the
equivalence relation \( R \). Denote by \( u_{1},\dots ,u_{n}\in W^{*}(X,R) \)
the canonical unitaries corresponding to these automorphisms. Then the numbers
\( \delta _{0}(R)=\delta _{0}(u_{1},\dots ,u_{n}) \), \( \delta ^{\omega }_{0}(R)=\delta _{0}^{\omega }(u_{1},\dots ,u_{n}) \)
and \( \delta ^{\omega }_{0,\kappa }(R)=\delta _{0,\kappa }^{\omega }(u_{1},\dots ,u_{n}) \)
depend only on \( R \).
\end{prop}
\begin{proof}
This follows from \ref{prop:independenceofgenerators}.
\end{proof}
In particular, \( \delta _{0}(R) \) is an invariant for the pair \( L^{\infty }(X)\subset W^{*}(X,R) \).

A more general statement holds: 

\begin{prop}
Let \( \Gamma  \) be a measurable \( r \)-discrete groupoid with base \( X \).
Assume that there exist a family of bisections \( \alpha _{1},\dots ,\alpha _{n} \),
which ``generates'' \( \Gamma  \) (i.e., so that every element of \( \Gamma  \)
can be written as the value of some product of \( \alpha _{1},\dots ,\alpha _{n},\alpha _{1}^{-1},\dots ,\alpha _{n}^{-1} \)).
Let \( u_{1},\dots ,u_{n}\in W^{*}(\Gamma ) \) be the unitaries canonically
associated to this family of bisections. Then the value of
\[
\delta _{0}(\Gamma )=\delta _{0}(u_{1},\dots ,u_{n})\]
depends only on \( \Gamma  \). The same statement holds true for \( \delta _{0}^{\omega } \)
and \( \delta ^{\omega }_{0,\kappa } \).
\begin{prop}
\label{prop:actionofz}Let \( \alpha  \) be a free action of \( \mathbb {Z} \)
on \( [0,1] \), which preserves Lebesgue measure \( \lambda  \). Denote by
\( A_{n}\subset [0,1] \) the set of points which have period exactly \( n \)
under the action \( \alpha  \); i.e., \( p\in A_{n} \) iff the set \( \bigcup _{k\in \mathbb {Z}}\alpha ^{k}(p) \)
has exactly \( n \) points. Denote by \( A_{\infty } \) the set \( [0,1]\setminus \bigcup _{n\geq 1}A_{n} \).
Let \( B=L^{\infty }[0,1] \) and \( M=B\rtimes _{\alpha }\mathbb {Z} \). Denote
by \( u \) the canonical unitary \( u\in M \), implementing \( \alpha  \).
Then
\[
\delta (u\btw B)=\sum _{n\geq 1}\frac{n-1}{n}\lambda (A_{n})+\lambda (A_{\infty }).\]
Moreover, if \( R_{\alpha } \) is the equivalence relation on \( [0,1] \)
induced by \( \alpha  \), then
\[
\delta (u\btw B)=C(R_{\alpha }),\]
where \( C(R_{\alpha }) \) is the cost of \( R_{\alpha } \) in the sense of
Gaboriau. The same conclusion holds for \( \delta  \) replaced by \( \delta ^{\omega } \),
\( \delta ^{\omega }_{\kappa } \), \( \delta _{0} \), \( \delta _{0}^{\omega } \)
and \( \delta _{0,\kappa }^{\omega } \).
\end{prop}
\end{prop}
\begin{proof}
Denote by \( p_{j}\in B \) the characteristic function of \( A_{j} \). Then
\( p_{j} \) commutes with \( u \) and \( \tau (p_{j})=\lambda (A_{j}) \).
By Proposition \ref{prop:additivityundercompression}, we obtain that
\[
\delta (u\btw B)=\sum _{j\geq 1}\lambda (A_{j})\delta (p_{j}up_{j}\btw p_{j}B)+\lambda (A_{\infty })\delta (p_{\infty }up_{\infty }\btw p_{\infty }B).\]
Note that \( p_{\infty }up_{\infty } \) has the same \( p_{\infty }B \)-valued
distribution as \( w(t)p_{\infty }up_{\infty } \), where \( w(t) \) is independent
from \( p_{\infty }B \) and free from \( p_{\infty }up_{\infty } \) with amalgamation
over \( p_{\infty }B \), and has the same \( * \)-distribution as free Brownian
motion started at identity and evaluated at time \( t \). It follows from Proposition
\ref{prop:stupidfdimbounds} that \( \delta (p_{\infty }up_{\infty }\btw B)\leq 1 \),
and from Proposition \ref{prop:approxexists} and Corollary \ref{corr:existanceofpermutation}
that \( \delta (p_{\infty }up_{\infty }\btw B)\geq 1 \). Hence \( \lambda (A_{\infty })\delta (p_{\infty }up_{\infty }\btw p_{\infty }B)=\lambda (A_{\infty }) \).

By Proposition \ref{prop:fdimcyclic}, we get that \( \delta (p_{j}up_{j}\btw p_{j}B)=\frac{j-1}{j} \)
for \( j<\infty  \). Hence \( \delta (u\btw B)=\sum _{n\geq 1}\frac{n-1}{n}\lambda (A_{n})+\lambda (A_{\infty }) \).
This is the same as the cost of \( R_{\alpha } \), see \cite{gaboriau:cost}.
\end{proof}
\begin{prop}
\label{prop:treeableequivrel}Let \( R \) be a treeable measurable equivalence
relation on \( [0,1] \), so that \( R=*_{i}R_{i} \), where each \( R_{i} \)
is generated by a single automorphism \( \alpha _{i} \). Let \( u_{j} \) be
the unitary in \( W^{*}(L^{\infty }[0,1],R) \), implementing \( \alpha _{j} \).
Then
\[
C(R)=\lim _{n\to \infty }\delta (u_{1},\dots ,u_{n}\btw B).\]

\end{prop}
\begin{proof}
We have
\[
C(R)=\sum _{i}C(R_{i})=\sum \delta (u_{i}\btw B)=\lim _{n\to \infty }\delta (u_{1},\dots ,u_{n}\btw B),\]
since \( u_{j} \) are free with amalgamation over \( B \).
\end{proof}
\begin{prop}
Let \( R \) be an equivalence relation possessing a finite graphing. Write
\( \delta (R)=\delta _{0,\kappa }^{\omega }(R) \). Then we have:
\begin{enumerate}
\item If \( R \) is the free product of equivalence relations \( R_{1},R_{2} \),
each having a finite graphing, then \( \delta (R_{1}*R_{2})=\delta (R_{1})+\delta (R_{2}) \)
\item If \( R \) is treeable, then \( \delta (R)=C(R) \), the cost of \( R \).
\item In general, \( \delta (R)\leq C(R) \).
\end{enumerate}
\end{prop}
\begin{proof}
The first and second properties follows from the additivity of \( \delta ^{\omega }_{0,\kappa } \)
for families of unitaries which are \( * \)-free over \( B \), and from Proposition
\ref{prop:treeableequivrel}. The last property follows from the fact that for
any finite graphing \( \alpha _{1},\dots ,\alpha _{m} \), denoting by \( R_{j} \)
the equivalence relation generated by \( \alpha _{j} \) and by \( u_{j} \)
the canonical unitary implementing \( \alpha _{j} \), we have:
\[
\sum C(R_{j})=\sum \delta _{0,\kappa }^{\omega }(u_{j})\geq \delta _{0,k}^{\omega }(u_{1},\dots ,u_{m})=\delta (R).\]
Choose now a measure-preserving automorphism \( \alpha  \) of \( [0,1] \),
so that \( \alpha  \) implements a free and ergodic action of \( \mathbb {Z} \),
and the induced equivalence relation \( R_{\alpha } \) is free from \( R \).
Let \( \bar{R}=R_{\alpha }\vee R \). Then
\[
\delta (\bar{R})=\delta (R_{\alpha })+\delta (R)=1+\delta (R)\leq \sum C(R_{\alpha _{j}})\]
 for any finite graphing \( \alpha _{1},\dots ,\alpha _{n} \) of \( \bar{R} \).
Let \( \beta _{1},\dots ,\beta _{n},\dots  \) be a graphing of \( R \). Then
\( \alpha ,\beta _{1},\dots ,\beta _{n},\dots  \) is a graphing of \( \bar{R} \).
If \( \sum C(R_{\beta _{j}})<+\infty  \), then there exists a finite graphing
\( \alpha ,\gamma _{1},\dots ,\gamma _{m} \) of \( \bar{R} \), with the same
cost as \( \alpha ,\beta _{1},\beta _{2},\dots  \) (indeed, given \( \beta _{i_{1}},\dots ,\beta _{i_{m}},\dots  \)
so that \( \sum \lambda (\textrm{domain}(\beta _{i_{k}}))\leq 1 \), one can
find integers \( n_{j} \) and \( m_{j} \) so that the domains and ranges of
\( \alpha ^{n_{j}}\beta _{i_{j}}\alpha ^{m_{j}} \), \( j=1,2,\dots  \) are
disjoint, and hence replace \( \beta _{i_{1}},\dots ,\beta _{i_{m}},\dots  \)
by a single automorphism, keeping the cost of the graphing the same). It follows
that
\[
\delta (R)=\delta (\bar{R})-1\leq C(R_{\alpha })+\sum _{j=1}^{m}C(R_{\gamma _{j}})-1=1+\sum _{j=1}^{\infty }C(R_{\alpha _{j}})-1,\]
since \( C(R_{\alpha })=1 \). Hence \( \displaystyle \delta (R)\leq \inf _{\alpha _{1},\dots ,\alpha _{m},\dots \, \textrm{graphing of }R}\sum _{j=1}^{\infty }C(R_{\alpha _{j}})=C(R) \).
\end{proof}

\subsection{Infinite number of generators.}

It is tempting to define, for an finite or infinite set \( S \) of unitaries
\( u_{1},u_{2},\dots \in \mathcal{N}(B) \) the quantity
\[
\underline{\delta }(S\btw B)=\lim _{k\to \infty }\delta _{0}(u_{1},u_{2},\dots ,u_{k}:u_{1},u_{2},\dots ,u_{k},u_{k+1},\dots \btw B).\]
(here set \( u_{k}=1 \) if \( k>|S| \)). By Proposition \ref{prop:fdiminpresence},
it follows that \( \underline{\delta }(S\btw B)=\delta (u_{1},\dots ,u_{n}\btw B) \)
if \( S=\{u_{1},\dots ,u_{n}\} \) is finite. In general, clearly \( \underline{\delta }\leq \delta  \).
We could not prove that \( \bar{\delta }(S\btw B) \) depends on the elements
of \( S \) only up to orbit-equivalence. In the case that \( u_{1},\dots ,u_{n},\dots  \)
form an infinite family, but are free with amalgamation over \( B \), one has

\begin{prop}
If \( u_{1},u_{2},\dots  \) are free with amalgamation over \( B \), then
\[
\underline{\delta }(\{u_{1},u_{2},\dots \}\btw B)=\lim _{n}\delta _{0}(u_{1},\dots ,u_{n}\btw B).\]
 
\end{prop}
\begin{proof}
We have that for each \( n \),
\[
\delta _{0}(u_{1},\dots ,u_{n}:u_{1},u,\dots \btw B)=\delta _{0}(u_{1},\dots ,u_{n}:u_{1},\dots ,u_{n}\btw B),\]
because \( u_{n+1},u_{n+2},\dots  \) are free from \( u_{1},\dots ,u_{n} \)
with amalgamation over \( B \). The rest follows from Proposition \ref{prop:fdiminpresence}.
\end{proof}

\section{Dynamical free entropy dimension of automorphisms.}

Let \( R \) be an equivalence relation on a measure space \( X \). We say
that \( \alpha  \) is an automorphism of \( R \), if \( \alpha  \) is an
automorphism of the von Neumann algebra \( W^{*}(X,R) \), so that \( \alpha (f),\alpha ^{-1}(f)\in L^{\infty }(X) \)
for all \( f\in L^{\infty }(X)\subset W^{*}(X,R) \). More generally, if \( M \)
is a von Neumann algebra, and \( B\subset M \) is a diffuse abelian subalgebra,
we say that \( \alpha  \) is an automorphism of \( B\subset M \), if \( \alpha (B)=B \).

For general automorphisms of a II\( _{1} \) factor \( M \), Voiculescu defined
its dynamical free entropy dimension in \cite[Section 7.2]{dvv:entropy2}. Unfortunately,
we don't at present know enough about free dimension \( \delta  \) to be able
to compute this invariant of an automorphism in all but very simple cases (for
example, it is trivial for any automorphism of the hyperfinite II\( _{1} \)
factor). It is natural, in view of relatively good behavior of \( \delta (\cdots \btw B) \)
with respect to orbit-equivalence operations, to try to use Voiculescu's definition
for automorphisms of groupoids, with the obvious modification of replacing \( \delta  \)
with \( \delta (\cdots \btw B) \). 

It will be useful to introduce the following notation. If \( F=(u_{1},\dots ,u_{n}) \),
then \( \alpha (F)=(\alpha (u_{1}),\dots ,\alpha (u_{n})) \), \( F_{k}=\cup _{j=0}^{k-1}\alpha ^{k}(F) \)
and \( F_{\infty }=\cup _{n=-\infty }^{\infty }\alpha ^{n}(F) \). For an automorphism
\( \alpha  \) of \( B\subset M \), set
\[
\underline{\delta }(\alpha ;F)=\limsup _{m}\frac{1}{m}\delta _{0}(F_{m}:F_{\infty }\btw B),\]

\[
\delta (\alpha ;F)=\limsup _{m}\frac{1}{m}\delta _{0}(F_{m}\btw B).\]

Note that \( \underline{\delta }<\delta  \).

\begin{defn}
Let \( \alpha  \) be an automorphism of \( B\subset M \). Define its dynamical
free entropy dimension to be
\[
\delta (\alpha )=\limsup _{F}\delta (\alpha ;F),\]
where \( F \) ranges over the set of all finite families of unitaries in \( \mathcal{N}(B) \),
ordered by inclusion.
\end{defn}
Since \( \delta _{0}(F\btw B)\leq |F| \), it follows that \( \delta (\alpha ;F)\leq |F| \).

\begin{defn}
Let \( F=(u_{1},\dots ,u_{n})\in \mathcal{N}(B)^{n} \). We say that \( F \)
is a \emph{weak generator for \( \alpha  \)}, if \( M=W^{*}(F_{\infty },B). \) 
\end{defn}
\begin{prop}
\label{prop:weakgen}Let \( \alpha  \) be an automorphism of an equivalence
relation \( R \), \( M=W^{*}(X,R) \) and \( B=L^{\infty }(X)\subset M \).
Let \( F \) be a weak generator, and let \( G\supset F \) be a finite family
of unitaries in the normalizer of \( B \). Then
\[
\delta (\alpha ;G)\leq \delta (\alpha ,F).\]

\end{prop}
\begin{proof}
Note that if \( F \) is a generator and \( F\subset F' \), then \( F' \)
is also a generator. It is thus sufficient to prove the Proposition for the
case that \( G=F\cup \{w\} \) for a single unitary \( w\in \mathcal{N}(B) \). 

If we replace \( F \) with \( F'=\alpha ^{-k}(F) \), then \( \delta (\alpha ;F)=\delta (\alpha ;F') \),
since \( \delta (u_{1},\dots ,u_{n}\btw B)=\delta (\alpha (u_{1}),\dots ,\alpha (u_{n})\btw B) \).

Let \( \rho >0 \) be fixed. Then there exists a projection \( p\in B \), \( N,M>0 \)
and a unitary \( v\in W^{*}(B,F_{N})\cap \mathcal{N}(B) \), so that \( pv\alpha ^{M}(w)^{*}=v\alpha ^{M}(w)^{*}p=p \)
and \( \tau (p)\geq 1-\rho  \). Write \( r=v\alpha ^{M}(w)^{*} \). Then
\begin{eqnarray*}
\delta _{0}(F_{m},w,\dots ,\alpha ^{m}(w)) & = & \delta _{0}(F_{m},w,\dots ,\alpha ^{m-1}(w)\btw B)\\
 & \leq  & \delta _{0}(F_{m},\alpha ^{M}(w),\dots ,\alpha ^{m-N-M-1}(w)\btw B)\\
 &  & +\delta _{0}(w,\dots ,\alpha ^{M-1}(w),\alpha ^{m-N}(w),\dots ,\alpha ^{m-1}(w)\btw B)\\
 & \leq  & \delta _{0}(F_{m},r,\dots ,\alpha ^{m-N-1-M}(r)\btw B)+N+M\\
 & \leq  & \delta _{0}(F_{m}\btw B)+\delta _{0}(r,\dots ,\alpha ^{m-N-M-1}(r)\btw B)+M+N\\
 & \leq  & \delta _{0}(F_{m}\btw B)+(m-N-M)\delta _{0}(r\btw B)+N+M.
\end{eqnarray*}
Since \( \delta _{0}(r\btw B)=\tau (1-p)\delta _{0}((1-p)r\btw B)+\tau (p)\delta _{0}(p\btw B)\leq \tau (1-p)\leq \rho  \),
we get
\[
\delta _{0}(\alpha ;F,w)\leq \delta _{0}(\alpha ;F)+\limsup _{m}\frac{1}{m}[(m-N-M)\rho +N+M]=\delta _{0}(\alpha ;F)+\rho .\]
Since \( \rho >0 \) is arbitrary, the conclusion follows.
\end{proof}
\begin{prop}
\label{prop:gen}Let \( \alpha  \) be an automorphism of \( B\subset M \).
Let \( F \) be a weak generator, and let \( G\supset F \) be a finite family
of unitaries in the normalizer of \( B \). Then \( \delta (\alpha ;G)\geq \underline{\delta }(\alpha ,F). \)
\end{prop}
\begin{proof}
Let \( H \) be a family so that \( G=F\cup H \). Then \( \delta _{0}(G_{m}\btw B)\geq \delta _{0}(F_{m},H_{m}:F_{\infty }\btw B)\geq \delta _{0}(F_{m}:F_{\infty }\btw B) \)
because \( H_{m}\subset W^{*}(B,F_{\infty }) \), so that Proposition \ref{prop:indepofgenineq2}
applies. Thus, by definition of \( \underline{\delta } \), we get that \( \delta _{0}(\alpha ;G)\geq \underline{\delta }(\alpha ;F). \)
\end{proof}
\begin{defn}
We say that a family \( F \) of unitaries in \( \mathcal{N}(B) \) is a \emph{generator
for \( \alpha  \)}, if it is a weak generator, and in addition
\[
\underline{\delta }(\alpha ;F)=\delta (\alpha ;F).\]
for all \( m\geq 0 \).
\end{defn}
\begin{prop}
If \( F \) is a generator for an automorphism of an equivalence relation \( R \),
then \( \delta (\alpha )=\delta (\alpha ;F) \).
\end{prop}
\begin{proof}
Using the fact that \( F \) is a generator and Proposition \ref{prop:gen},
we find that for any family \( G\supset F \), \( \delta (\alpha ;G)\geq \underline{\delta }(\alpha ,F)=\delta (\alpha ,F) \).
Combining this with Proposition \ref{prop:weakgen} gives \( \delta (\alpha ;F)\geq \delta (\alpha ;G)\geq \delta (\alpha ;F) \).
It follows that \( \delta (\alpha ;G)=\delta (\alpha ;F) \). It follows that
\( \limsup _{H}\delta (\alpha ;H)=\delta (\alpha ;F) \) since \( F\subset H \)
for sufficiently large \( H \).
\end{proof}

\subsection{Examples of automorphisms.}

We conclude by giving an example for which the dynamical free entropy dimension
invariant is non-trivial. Let \( \alpha  \) be a free measure-preserving action
of the free group \( \mathbb {F}_{n} \) on a finite measure-space \( X \)
(e.g., one can take the Bernoulli action of \( \mathbb {F}_{n} \) on \( \prod _{g\in \mathbb {F}_{n}}\{0,1\} \)).
Let \( Q \) be the associated equivalence relation. Denote by \( g_{1},\dots ,g_{n} \)
be infinite free generators of \( \mathbb {F}_{n} \). Consider the equivalence
relation \( R \) induced on \( X \) by the action of the subgroup \( G \)
of \( \mathbb {F}_{n} \) generated by the set \( \{g_{n}^{k}g_{j}g_{n}^{-k}:1\leq j\leq n-1,\, k\in \mathbb {Z}\} \).
It is not hard to see that \( G\cong \mathbb {F}_{\infty } \). Then \( W^{*}(X,R)\subset W^{*}(X,Q) \).
Denote by \( w\in W^{*}(X,Q) \) the unitary implementing the action of \( g_{n} \).
Then \( wW^{*}(X,R)w^{*}=W^{*}(X,R) \) and \( wL^{\infty }(X)w^{*}=L^{\infty }(X) \).
It follows that \( \alpha (y)=wyw^{*} \) is an automorphism of \( W^{*}(X,R) \),
and moreover is an automorphism of the equivalence relation \( R \). This automorphism
is called a free shift of multiplicity \( n-1 \).

Denote by \( u_{i}\in W^{*}(X,R) \) the unitary implementing the action of
\( g_{i} \), \( 1\leq i\leq n-1 \). Set \( F=(u_{1},\dots ,u_{n-1}) \).

\begin{claim}
\( F \) is a generator for \( \alpha  \).
\end{claim}
\begin{proof}
Note that \( \alpha ^{k}(u_{j}) \) is the unitary corresponding to \( g_{n}^{k}u_{j}g_{n}^{-k} \).
It follows that \( F_{\infty } \) together with \( B \) generates \( W^{*}(X,R) \).
Hence \( F \) is a weak generator.

For any fixed \( m \), we have \( \delta _{0}(F_{m}:F_{\infty }\btw B)=\delta _{0}(F_{m}:F_{m}\btw B) \)
since the set \( \{\alpha ^{k}(u_{j}):1\leq j\leq n-1,\, k\notin \{0,\dots ,m-1\}\} \)
is free from \( F_{m} \) with amalgamation over \( B \) (see Proposition \ref{prop:inpresencefree}).
On the other hand, \( \delta _{0}(F_{m}:F_{m}\btw B)=\delta _{0}(F_{m}\btw B) \)
by Proposition \ref{prop:fdiminpresence}. Hence \( \delta _{0}(F_{m}:F_{\infty }\btw B)=\delta _{0}(F_{m}\btw B) \),
and thus \( \delta (\alpha ;F)=\underline{\delta }(\alpha ;F) \). Hence \( F \)
is a generator.
\end{proof}
\begin{claim}
\( \delta (\alpha )=n-1 \).
\end{claim}
\begin{proof}
We have that \( \delta (\alpha )=\delta (\alpha ;F) \), since \( F \) is a
weak generator. But \( \delta _{0}(F_{m}\btw B)=\sum _{i=1}^{n-1}\sum _{k=0}^{m-1}\delta (\alpha ^{k}(u_{j})\btw B) \),
since \( \{\alpha ^{k}(u_{j})\}_{j,k} \) are free with amalgamation over \( B \).
Since each \( \alpha ^{k}(u_{j}) \) implements a free action of the integers,
\( \delta _{0}(\alpha ^{k}(u_{j})\btw B)=1 \), so that \( \delta _{0}(F_{m}\btw B)=m(n-1) \).
Thus \( \delta (\alpha ;F)=n-1 \).
\end{proof}
By replacing in the construction above the set \( X \) by the set \( Z=X\sqcup Y \),
so that \( \mu _{Z}(X)=t \), \( \mu _{Z}(Y)=1-t \), and letting \( \mathbb {F_{n}} \)
act trivially on \( Y \), one obtains examples of automorphisms of equivalence
relation having dynamical free entropy dimension \( tn \). By varying \( t \),
it is clear that one can obtain all numbers in \( (0,+\infty ) \) as values
of dynamical free entropy dimension of an automorphism of an equivalence relation.
Clearly, \( \delta (\textrm{id})=0 \), so in fact all numbers in \( [0,+\infty ) \)
can be obtained. We don't know if the infinite-multiplicity free shift has dynamical
free entropy dimension \( +\infty  \), although we suspect this is the case.

\subsection{Groups.}

It is possible to define an invariant for group automorphisms in the same way.
If \( G \) is a group and \( g\in G \), denote by \( u(g) \) the unitary
in the group von Neumann algebra of \( G \), corresponding to \( g \). Let
\( \alpha  \) be an automorphism of \( G \). For a finite family \( F=(g_{1},\dots ,g_{n}) \)
in \( G \), define \( \alpha (F) \), \( F_{m} \) and \( F_{\infty } \) in
the obvious way. Set \( \delta (\alpha ;F)=\limsup _{m}\frac{1}{m}\delta _{0}(F_{m}) \)
and \( \underline{\delta }(\alpha ;F)=\limsup _{m}\frac{1}{m}\delta _{0}(F_{m}:F_{\infty }) \)
(here \( \delta _{0} \) is the modified free entropy dimension of Voiculescu,
see \cite{dvv:improvedrandom}). The results of this section, after an appropriate
modification, remain true in this case. We leave the details to the reader,
but summarize the results.

\begin{enumerate}
\item Say that a family \( F \) of elements of \( G \) is a weak generator for \( \alpha  \),
if \( F_{\infty } \) generates \( G \). Say that \( F \) is a generator for
\( \alpha  \), if it is a weak generator of \( G \), and in addition \( \underline{\delta }(\alpha ;F)=\delta (\alpha ;F). \)
\item If \( F \) is a generator, then \( \delta (\alpha )=\delta (\alpha ;F) \).
\item Let \( H \) be a group and \( G=*_{i\in \mathbb {Z}}H \). Let \( \alpha  \)
be the free shift automorphism. Then \( \delta (\alpha )=\delta (H) \).
\end{enumerate}
More generally, the results of this section remain valid for automorphisms of
\( r \)-discrete finite measure groupoids. We leave the details to the reader.

\bibliographystyle{amsplain}

\providecommand{\bysame}{\leavevmode\hbox to3em{\hrulefill}\thinspace}

\end{document}